\documentclass[11pt]{amsart}
\usepackage[all]{xy}
\usepackage{mathrsfs,amssymb}

\theoremstyle{plain}
\newtheorem{theorem}{Theorem}[section]
\newtheorem{proposition}[theorem]{Proposition}
\newtheorem{lemma}[theorem]{Lemma}
\newtheorem{corollary}[theorem]{Corollary}

\theoremstyle{definition}
\newtheorem{definition}[theorem]{Definition}
\newtheorem{example}[theorem]{Example}
\newtheorem*{acknowledgement}{Acknowledgements}
\newtheorem{remark}[theorem]{Remark}

\theoremstyle{remark}

\allowdisplaybreaks

\newcommand\defn[1]{{\bf #1}}
\newcommand\simp{{}^R_{\rm sim}}
\newcommand\virt{^{\rm virt}}

\newcommand{\wt}{\widetilde}
\newcommand{\inv}{^{-1}}
\newcommand{\vb}{(V,\beta)}

\newcommand{\Bl}{\operatorname{Bl}}
\newcommand{\mult}{\operatorname{mult}}
\newcommand{\ra}{\rightarrow}

\newcommand{\rarr}{\rightarrow}

\newcommand{\PGL}{\operatorname{PGL}}

\newcommand{\Pic}{\operatorname{Pic}}
\newcommand{\Spec}{\operatorname{Spec}}
\newcommand{\Proj}{\operatorname{Proj}}

\newcommand{\Aut}{\operatorname{Aut}}
\newcommand{\Cone}{\operatorname{Cone}}

\newcommand{\bC}{{\mathbb{C}}}
\newcommand{\bG}{{\mathbb{G}}}

\newcommand{\bN}{{\mathbb{N}}}
\newcommand{\bP}{{\mathbb{P}}}
\newcommand{\bA}{{\mathbb{A}}}
\newcommand{\bQ}{{\mathbb{Q}}}
\newcommand{\bR}{{\mathbb{R}}}
\newcommand{\bZ}{{\mathbb{Z}}}

\newcommand{\cB}{{\mathcal{B}}}
\newcommand{\cC}{{\mathcal{C}}}
\newcommand{\cD}{{\mathcal{D}}}

\newcommand{\cI}{{\mathcal{I}}}
\newcommand{\cL}{{\mathcal{L}}}

\newcommand{\cO}{{\mathcal{O}}}

\newcommand{\cW}{{\mathcal{W}}}

\newcommand{\oC}{\overline C}

\newcommand{\oM}{\overline M}

\newcommand{\des}[3]{\left\langle{#1}\right\rangle^{#2}_{#3}}

\newcommand{\sbul}{s_{\bullet}}

\newcommand{\bt}{{\boldsymbol{t}}}
\newcommand{\bk}{{\boldsymbol{k}}}

\newcommand{\Ap}{A^{+}}
\newcommand{\Am}{A}

\newcommand{\Ds}{D_{\sigma}}

\title[Weighted gravitational descendants]%
{Moduli of weighted stable maps and their gravitational descendants}
\date{September 18, 2006; revised November 2, 2007}
\author{Valery Alexeev}
\address{Department of Mathematics, University of Georgia,
  Athens GA 30602, USA}
\email{valery@math.uga.edu}
\author{G. Michael Guy}
\address{Department of Mathematics, University of Georgia, Athens GA
  30602, USA} \email{{guy@math.uga.edu}}

\begin{document}
\maketitle
\begin{abstract}
  We study the intersection theory on the moduli spaces of maps of
  $n$-pointed curves $f:(C,s_1,\dots s_n)\to V$ which are stable with
  respect to the weight data $(a_1,\dotsc, a_n)$, $0\le a_i\le 1$.
  After describing the structure of these moduli spaces, we prove a
  formula describing the way descendant invariants change under a wall
  crossing. As a corollary, we compute
  the weighted descendants in terms of the
  usual ones, i.e. for the weight data $(1,\dotsc,1)$, and vice versa.
\end{abstract}

\tableofcontents

 \section*{Introduction}
 \label{sec:intro}

Moduli spaces $\oM_{g,n}(V,\beta)$ of stable maps $f:(C,s_1,\dotsc,
s_n)\to V$ from $n$-pointed curves to an algebraic variety have been
extensively studied. In particular, they are used to define the
Gromov-Witten invariants of $V$ and quantum cohomology. In the
absolute case, when $V$ is a point, these are the
Deligne-Mumford-Knudsen's moduli spaces of stable $n$-pointed
curves.

B. Hassett \cite{Hassett_WeightedStableCurves} defined weighted
analogs of stable curves. In this version, one attaches to every
point $s_i$ a weight $0\le a_i\le 1$, and modifies the definition of
a stable curve accordingly. Hassett constructed the corresponding
moduli spaces $\oM_{g,A}$ and proved a number of their nice
properties: they are all Deligne-Mumford stacks, smooth if all
$a_i>0$, and their coarse moduli spaces are projective. For two
weight data $A\ge B$, there is a reduction morphism
$\rho_{B,A}:\oM_{g,A}\to \oM_{g,B}$; it is birational, and Hassett
described the exceptional divisors which it contracts.

It is natural to also define a weighted version of a stable map. In
fact, the definition can be given for higher dimensions as well--see
\cite{Alexeev_ICM} for a more detailed discussion on this subject
(the present paper answers \cite[Question 7.3]{Alexeev_ICM}). In the
unweighted case, the stable maps from surfaces were defined and
studied in \cite{Alexeev_Mgn}. Hacking studied surface pairs
$\big(X,(3/d+\epsilon)D\big)$ with weights in \cite{Hacking_PlaneCurves}, and
the moduli compactifications of toric and abelian pairs
\cite{Alexeev_CMAV} can be interpeted as moduli spaces of pairs
$(X,\Delta+\epsilon D)$, resp. $(X,\epsilon D)$ with weights.

In this paper, we first construct the moduli spaces $\oM_{g,A}(V,
\beta)$ of maps $f:(C,s_1,\dotsc, s_n)\to V$ from $n$-pointed curves
to a projective variety $V$ (or, more generally, a flat projective
scheme over the base scheme) with respect to the weight data
$A=(a_i)$, $0\le a_i\le 1$. We give two constructions.

The advantage of the first one that it is a very elementary reduction
to the usual, unweighted case established in \cite{BehrendManin}.
However, for fixed $V,\beta,n$, \cite{BehrendManin} only gives a
Deligne moduli stack 
outside of finitely many characteristics (``bounded by the
characteristic'' in the language of \cite{BehrendManin}).

The second construction is more general, and gives an Artin moduli
stack and a coarse moduli space over an arbitrary locally Notherian
scheme.  This construction, which follows
\cite{Mumford_GIT,Kollar_Projectivity}, was previously used for
surfaces, and can be used in higher dimensions as well (with many
significant technicalities, which we do not discuss here).

Next, we define the psi classes on $\oM_{g,A}(V,\beta)$. In the case
when all $a_i>0$, this is done in the usual way by restricting the
cotangent bundle to a section. Since the section passes through the
locus of the universal family where the morphism is smooth, the psi
classes are invertible sheaves on the moduli stack. In the case of
some $a_i=0$, we adjust this construction slightly.

In \cite{Hassett_WeightedStableCurves}, Hassett proved that the
parameter space of the possible weight data $A$, which for $g\ge2$
is simply $[0,1]^n$, is divided into finitely many ``chambers'',
within which the moduli spaces $\oM_{g,A}$ do not change. For each
fixed $(g,n,V,\beta)$, we define a similar
chamber decomposition for the spaces $\oM_{g,A}(V,\beta)$.

It is easy to see that one can move from one chamber to another by a
sequence of ``simple'', generic wall crossings. We organize this
data as follows. To each weight data, or chamber, we associate a
simplicial complex with vertices $\{1,\dotsc,n\}$. (For example, the
complex for the weight data $(1^n)$ is $n$ disjoint points, and for
the weight data $(0^n)$ it is the $(n-1)$-dimensional simplex.) Then
crossing a simple wall corresponds to adding or removing a single
simplex to or from this complex.

 Next, we restrict to a case when $V$ is a
smooth projective variety defined over a field. We define
gravitational descendants $\des{\tau_{k_1}(\gamma_1)\cdots
\tau_{k_n}(\gamma_n)}{}{g,A}$ of~$V$. We then prove a  formula
expressing the way each gravitational descendant changes under a
simple wall crossing. The formula says that the difference is again
a gravitational descendant, but on a simpler space with fewer marked
points. In particular, it also implies that the ordinary
Gromov-Witten invariants coincide with their weighted analogues. As a
corollary, we get formulas for $\des{\tau_{k_1}(\gamma_1)\dots
\tau_{k_n}(\gamma_n)}{}{g,A}$ in terms of the ``unweighted''
descendants, and vice versa.

Many authors have studied the change in gravitational descendants
under ``abrupt'' moves, when going from $n$-pointed curves to
$(n-1)$-pointed curves, etc. The Witten's formulas in \cite{Witten}
relating products of $\kappa$ classes and $\psi$ classes are one of
the first examples, and \cite{ArbarelloCornalba} contains many more.
Definitions were given for the ``canonical'' descendants that behave
well under the forgetful morphism, as the pullbacks from a space
with fewer marked points. For example, Graber-Kock-Pandharipande
\cite{GKP} introduced ``modified psi classes'' which can be
interpreted as pullbacks of psi classes for the weight data $(0^n)$
for $g>0$ or  $\beta\neq 0$, and this construction was extended to
the remaining case in \cite{Kock}.

In these terms, our formulas can be understood as describing the
change under the ``micro'' moves, and filling out the entire cube
$[0,1]^n$ rather than jumping along the edges. The previously studied
invariants can be understood as the invariants at the corners, i.e.
points with coordinates all $1$s and $0$s.  In particular, we prove
that the Morita-Miller-Mumford's kappa numbers are simply the weighted
descendants for the weight data $(0^n)$, or, equivalently,
$(\epsilon^n)$.

In the last section, we study the weighted analogs of the string,
dilaton, and divisor equations.

\medskip We note a related and independent work \cite{BayerManin},
which appeared on arXiv shortly before ours.  Arend Bayer
and Yuri~I. Manin also construct the moduli of weighted stable maps
and establish the basic properties of reduction morphisms (luckily,
the proofs are fairly different). Next, they study and prove some
axioms of the virtual fundamental class.  Then they study the way the
stability conditions in derived categories change under wall
crossings.

In contrast, our focus is on the weighted gravitational descendants,
and we consider the wall crossing formula in Theorem
\ref{theorem:wallcross} to be the central result of this paper.

\begin{acknowledgement}
%
%
  The authors thank Rahul Pandharipande for offering useful
  references. We are grateful to the referee for many very useful comments.
  This research was partially supported by NSF.

  The second author would also like to thank Jim Bryan who led a group
at the 2005 Summer Institute in Algebraic Geometry in which he was a
member. It is his notes which he distributed during this session
that have served as the starting point for much of the second
author's learning about descendants.

\end{acknowledgement}

\section{Moduli of weighted stable maps}
\label{sec:moduli}

For the definitions of individual varieties and pairs, we work over
an algebraically closed field. For the definitions of families and
moduli functors, we work over a Noetherian base scheme
$\cB$, and all products are fibered products over $\cB$. $V$ will
denote a flat projective scheme over $\cB$.

\begin{definition}
  A \defn{weight} in this paper will be a real number $0\le a_i\le 1$.
  A \defn{weight data} will be an ordered $n$-tuple
  $A=(a_1,\dotsc,a_n)$ of weights. We will call a weight data positive if each weight is positive.
\end{definition}

We will use abbreviations: for example, we will write $(1^n)$ for
the weight data $(1,\dotsc,1)$, and $(1^2,\epsilon^{n-2})$ for
$(1,1,\epsilon,\dotsc, \epsilon)$.

\begin{definition}\label{defn:stable-map}
  Over an algebraically closed field,
  a \defn{stable map for the weight data $\mathbf{A}$}, or an
  \defn{$\mathbf{A}$-stable map} is a proper morphism $f:C\to V$ from
  a connected reduced curve $C$ to a scheme $V$, together with $n$
  \emph{ordered} points $s_1,\dotsc, s_n\in C$ which satisfies the
  following two conditions:
  \begin{enumerate}
  \item (on singularities) $C$ has at most nodes; for every smooth
    point $P\in C$, the multiplicity
    \begin{displaymath}
      \mult_P = \sum_{s_i = P} a_i \le 1,
    \end{displaymath}
    and for a node $P\in C$, one has $\mult_P =0$;
  \item (numerical) the $\bR$-line bundle $\omega_C(\sum a_i s_i)$ is
    $f$-ample, i.e. for every irreducible component $E$ of $C$
    collapsed by $f$ to a point, one has
    \begin{displaymath}
      \deg \omega_C\left(\sum a_i s_i\right)\Big{|}_E =
      2p_a(E)-2 + |E\cap (C-E)| + \sum_{s_i\in E} a_i >0.
    \end{displaymath}
  \end{enumerate}
  A \defn{stable curve} is a stable map to a point.
\end{definition}

Thus, the points with $a_i=0$ can coincide with the nodes, but
points with positive weights cannot. As usual, the numerical
condition is only a restriction on the collapsed components $E$
which are $\bP^1$, elliptic curves, or rational curves with a single
node.

\begin{definition}\label{definition:fammaps}
  Let $V$ be a flat scheme over $\cB$.
  A family of stable maps over $S$ is a morphism of schemes
  $f:C\to V\times_{\cB} S$ together with sections $g_i:S\to C$ of
  $\pi=p_2\circ f:C\to S$, such that
  \begin{enumerate}
  \item $\pi:C\to S$ is flat,
  \item   every geometric fiber
  $\big(C, s_i=g_i(S)\big)_{\bar s} \to V_{\bar s}$
  is an $A$-stable map.
  \end{enumerate}
\end{definition}

\begin{definition}
  The moduli stack $\oM_{g,A}(V)$ associates to every scheme $S/\cB$
  a category whose objects are families of stable
  maps over $S$ such that every curve $C_{\bar s}$ has arithmetic
  genus $g$; and arrows are isomorphisms of families over $S$.

  The moduli functor is defined by associating to $S$ the set of such
  families modulo isomorphisms.
\end{definition}

\begin{remark}\label{rem:ratnl-weights}
  For any weight data $(a_i)$, one may decrease the nonzero $a_i$ ``a
  little" without changing the stability conditions: a family is stable
  for $(a_i)$ iff it is stable for $(a_i-\epsilon_i)$. Hence, for
  proving statements about families of $A$-stable maps, one can always
  assume that each $a_i$ is rational.
\end{remark}

Stable weighted curves were defined by Hassett
\cite{Hassett_WeightedStableCurves} who constructed their moduli
spaces and gave a detailed description.
The extension to the relative case, of course, is immediate (see for
example \cite{Alexeev_Mgn} where the moduli of stable maps \emph{from
  surface pairs} were constructed, although without weights).

For completeness, let us give the definition in the case when points
$s_i$ are replaced by ``something'' of higher degree:

\begin{definition}
  Fix positive integers called degrees $d_i$, $i=1,\dotsc,n$.  A
  family of stable maps over $S$ for the data $\big(A,(d_i)\big)$ is a
  morphism of schemes $f:C\to V\times S$ together with morphisms
  $g_i:D_i\to C$ such that
  \begin{enumerate}
  \item the compositions $\pi\circ g_i:D_i\to S$
  are finite and \'etale.
\item on every geometric fiber, the map $f_{\bar s}:C_{\bar s}\to
  V_{\bar s}$ is stable with respect to the $\sum d_i$ points
  $g_i(D_{i})_{\bar s}$ for the weight data in which each
  $a_i$ is repeated $d_i$ times.
  \end{enumerate}
\end{definition}

Thus, the only difference with the previous case is that we forget the
order in each group of $d_i$ points. It is obvious that the
corresponding moduli stack is the $S_{d_1}\times\dotsb\times
S_{d_n}$-quotient of the previous moduli stack, and the same holds for the
coarse moduli spaces.

\begin{remark}
  Note, however, that one gets a very different moduli stack if
  instead of finite maps $g_i:D_i\to C$ one considers divisors or
  subschemes $D_i\subset C$ of degree $d_i$. When $D_i$ lie in the
  smooth part of $C$, the coarse moduli space is the same, but the
  moduli stacks are different--ours have extra automorphisms. Over the
  nodes, the problems multiply, of course.
\end{remark}

We will now fix a projective scheme $V\subset \bP^N$ with a very
ample sheaf $\cO_V(1)$, an integer $g\ge0$ and the weight data $A$.
As in the unweighted case, we will further subdivide $\oM_{g,A}(V)$
into a disjoint union
\begin{displaymath}
  \oM_{g,A}(V) = \coprod_{\beta} \oM_{g,A}(V,\beta)
\end{displaymath}
with pieces $\oM_{g,A}(V,\beta)$ of finite type.
There are several meanings that can be given to $\beta$:
\begin{enumerate}
\item When $V$ is a complex variety, one can fix a homology class\linebreak
  $\beta\in H_2(V,\bZ)$.
\item When the base is a field and $V$ is a variety, the following
  algebraic analog of the homology class is commonly used (see
  e.g. \cite{BehrendManin}):
\begin{displaymath}
  \beta:\Pic V\to \bZ, \quad L\mapsto \deg f^*(L),
\end{displaymath}
This can be extended to the case when the relative Picard scheme of
$V/\cB$ is nice, for example when $V=\bP^N_{\cB}$.
\item But we also would like to work with the case when $V$ is just a
  flat projective scheme over $\cB$. In this case, the simplest
  invariant which is locally closed on the fibers of a family
  and guarantees finite type is
  $$ \deg \beta = \deg f^* (\cO_V(1))\ge0. $$
  By abuse of notation, we will still use
  $\oM_{g,A}(V,\beta)$ instead of the more accurate but cumbersome
  $\oM_{g,A}(V,\deg\beta)$.
\end{enumerate}
\begin{remark}
Whenever the omission of any of $g,A,V,\beta$ is unlikely to lead to confusion, we may do so.
\end{remark}

\begin{theorem}\label{thm:stack-exists}
  The moduli stack $\oM_{g,A}(V,\beta)$ is a proper algebraic Artin
  stack with finite stabilizer.  For each $(g,n,V,\beta)$, there
  exists $N\in \bN$ such that
  \linebreak$\oM_{g,A}(V,\beta)\times\bZ[1/N]$ (i.e. outside of the
  finitely many positive characteristics dividing $N$) is a
  Deligne-Mumford stack.
\end{theorem}

By \cite{KeelMori}, this implies:
\begin{corollary}
  $\oM_{g,A}(V,\beta)$ has a coarse moduli space, a proper algebraic
  space over the base scheme $\cB$.
\end{corollary}
As the referee pointed out to us, the formation of this coarse moduli
space may not commute with the change of the base scheme
$\cB'\to\cB$, cf. \cite{AbramovichOlssonVistoli}.

This theorem is a baby case of a construction which can be carried
out for maps $f:(X,D_i)\to V$ with higher-dimensional $X$.
We will give two proofs. The first proof will be a
reduction to the case of weight data $(1^n)$, i.e. the ``standard
unweighted'' case, where the statement is well-known, see
e.g. \cite{BehrendManin}, although in a weaker form, outside of
finitely many characteristics. The second proof is more general.

\begin{remark}
  Let $f:\bP^1\to \bP^1$ be the Frobenius map $x\mapsto x^p$, and let
  us work in characteristic $p>0$. Then $f$ is a stable map whose
  automorphism group scheme is $\mu_p=\Spec k[x]/(x^p-1)$ is not
  reduced. This shows that the stack $\oM_{0,0}(\bP^1,p\bP^1)$ is not
  Deligne-Mumford in general.
\end{remark}

\begin{theorem}\label{thm:properness}
  The stack $\oM_{g,A}(V,\beta)$ satisfies the valuative criterion of
  properness: every family over a punctured $\Spec$ of a DVR or a
  punctured regular curve $S\setminus 0$ has at most one extension, and
  the extension always exists after a finite base change $(S',0)\to
  (S,0)$.
\end{theorem}
\begin{proof}
The proof is exactly the same for the weighted or unweighted case,
maps or curves (and indeed for higher-dimensional pairs as well),
and is a variation of \cite[\S2]{DeligneMumford}, cf. also
\cite[Prop.3.7]{Hassett_WeightedStableCurves}.

Let $U=S\setminus 0$, and let $C_U\to V\times U$ be the family.
Assume first that the generic fiber $C_{\eta}$ is smooth. By the
semistable reduction theorem \cite{ArtinWinters} (which holds in
complete generality, including the case of mixed characteristics),
after a base change $S'\to S$, there exists a semistable extension
$C'\to S'$ such that: $C'/S$ is flat and proper, $C'$ is regular,
and the central fiber is a reduced divisor with normal crossings.
We then proceed with the log Minimal Model Program for the
pair $(C',\sum a_is'_i)$ relative over $V\times S$, which is available
in the necessary generality by \cite{Shafarevich_LecsMinModels,
Lichtenbaum_CurvesDVR}.

First, we inductively contract $f$-exceptional curves $E$ in the
central fiber with $(K_{C'}+\sum a_is'_i)E<0$, where $K_C'$ stands
for the dualizing sheaf $\omega_{C'}$, which is invertible. These
are easily seen to be $(-1)$-curves, i.e.  forms of $\bP^1$ with
$E^2=-1$. After that, we contract all the $f$-exceptional curves
with $(K_{C'}+\sum a_is'_i)E=0$. These are seen to be chains of
$(-2)$-curves not meeting the sections $s'_i$ with $a_i>0$. The
resulting surface $\oC'$ is Gorenstein and has $A_n$-singularities
of the form $xy=u t^{n+1}$, where $t$ is the uniformizing parameter
of the DVR, and $u$ is a unit.  Hence, the central fiber $\oC'_0$ is
a nodal curve, and $(\oC',s'_i)\to S'$ is the required extension.

Conversely, let $C\to S$ be an extension. From the description of
the deformations of ordinary double points, we see that $C$ has
singularities of the form $xy=ut^{n+1}$ and is Gorenstein. The
minimal resolution adds a chain of $(-2)$-curves. Let $\wt C$ be any
regular surface extending $C_U$ and dominating $C$, $F:\wt C\to C$.
Then $F$ is a sequence of simple blowups and the local singularity
condition on $C$ implies that
\begin{displaymath}
  K_{\wt C}+\sum a_i s_i = F^*\left(K_C+\sum a_is_i\right)+ \sum b_jE_j,
\end{displaymath}
where $E_j$ are $F$-exceptional and $b_j\ge 0$. This implies that
\begin{displaymath}
  C = \Proj_{V\times S} \bigoplus_{d\ge 0}
  \pi_*\cO_{\wt C} \big( d (K_{\wt C}+\sum a_i s_i) \big).
\end{displaymath}
In this formula, we assume that $a_i$ are rational by
Remark~\ref{rem:ratnl-weights}, and $d\in N\bN$ such that $Na_i\in
\bZ$.

Since any two extensions of $C_U$ can be dominated by a common model
$\wt C$, this implies the uniqueness of the extension.

The general case of a non-normal generic fiber is reduced to the
previous case by taking the normalization. After a finite base change, the normalization is a
disjoint union of the families as above with extra sections of weight
1. Given the extension for the
normalization, the extension for the original family is obtained in a
unique way by gluing along these sections.
\end{proof}

\begin{lemma}\label{lem:weights0and1}
  Let $A=(a_1,\dotsc,a_n)$ be a weight data and let $A\cup
  0^m=(a_1,\dotsc,a_n,0,\dotsc,0)$ be the weight data obtained by adding
  $m$ zeros.  Assume that Theorem~\ref{thm:stack-exists} has been
  proved for $A$, and let $(\cC,s_i)\to \oM_{g,A}(V,\beta)$ be its
  universal family.  Then Theorem~\ref{thm:stack-exists} holds for
  $A\cup 0^m$, and
  \begin{displaymath}
    \oM_{g,A\cup 0^m}(V,\beta) = \cC^m_{\oM_{g,A}(V,\beta)}
  \end{displaymath}
  is the $m$-th fibered power of $C$ over $\oM_{g,A}(V,\beta)$.
\end{lemma}
\begin{proof}
  Indeed, the only difference between families in $\oM_{g,A}(V,\beta)$
  and in $\oM_{g,A\cup 0^m}(V,\beta)$ is $m$ \emph{arbitrary} sections.
\end{proof}

Thus, we can and will \emph{assume that the weight $A$ is positive} until
the end of the first proof of Theorem~\ref{thm:stack-exists}.

\medskip
Each nodal $n$-pointed curve defines an equivalence relation $R(C,s_i)$
on the set $\{1,\dotsc, n\}$: $i\sim j$ iff $s_i=s_j$.

\begin{definition}\label{defn:R-compatible}
Let $R$ be an equivalence relation on $\{1,\dots,n\}$.
A map $f:(C,s_i)\to V$ with nodal~$C$ is called
\defn{$R$-compatible} if
  \begin{enumerate}
\item If $i\sim j$, then $s_i,s_j$ lie on the same
  irreducible component of $C$; \label{defn:R-compatible1}
\item $s_i$ are not nodes; \label{defn:R-compatible2}
\item if $i\not\sim j$, then $s_i\ne s_j$. \label{defn:R-compatible3}
\end{enumerate}
\end{definition}

\begin{lemma}
  The stack $U_{g,A}(V,\beta,R)$ of $A$-stable families such that
  every geometric fiber is $R$-compatible is an open substack of
  $\oM_{g,A}(V,\beta)$.
\end{lemma}
\begin{proof}
  Indeed, it is obvious that all three parts of Definition
  \ref{defn:R-compatible} are open
  conditions in families of nodal curves.
\end{proof}

\begin{definition}
  We define the \defn{$R$-simplified weight data $A\simp$} as follows:
  for each group of $R$-equivalent indices $i$, we set one of the
  weights $a_i$, for example for the smallest $i$, to $1$ and the
  other weights in that group to~$0$.
\end{definition}

\begin{lemma}\label{lem:simplify-weights}
  Let $f:(C,s_i)\to V$ be an $R$-compatible map with nodal $C$.
  Assume that for every $R$-equivalence class $r$ one has $\sum_{i\in
    r} a_i\le 1$.  Then if $f$ is $A$-stable, then it is
  $A\simp$-stable.
\end{lemma}
\begin{proof}
  By the conditon $\sum_{i\in r}a_i\le 1$ and the requirement
  \ref{defn:R-compatible}\eqref{defn:R-compatible3}, $(C,s_i)$
  satisfies the singularity condition for stability.

  Requirement \ref{defn:R-compatible}\eqref{defn:R-compatible1}
  implies that the numerical condition for weight $A\simp$ is the
  same, or better, than for weight $A$.
\end{proof}

\begin{corollary}\label{cor:local-structure}
  The stacks $\oM_{g,A}(V,\beta)$ and $\oM_{g,A\simp}(V,\beta)$ are
  locally isomorphic. Namely, for each point $(C,s_i,f)\in
  \oM_{g,A}(V,\beta)$ there exists an open neighborhood isomorphic to
  an open substack in $\oM_{g,A\simp}(V,\beta)$, where $R=R(C,s_i)$.
\end{corollary}
\begin{proof}
  For every $A$-stable map $f$, we take $R=R(C,s_i)$. Then we have an
  open neighborhood $U_{g,A}(V,\beta,R)\subset \oM_{g,A}(V,\beta)$,
  and by Lemma~\ref{lem:simplify-weights}, $U_{g,A}(V,\beta,R)$ is an
  open substack of $U_{g,A\simp}(V,\beta,R)\subset
  \oM_{g,A\simp}(V,\beta)$.
\end{proof}

\begin{proof}[First proof of Theorem~\ref{thm:stack-exists}]
  For the weight data $(1^n)$, i.e. in the standard unweighted case,
  Theorem~\ref{thm:stack-exists} is well known, see
  e.g. \cite{BehrendManin}. For the weight data
  $A\simp$, it follows by Lemma~\ref{lem:weights0and1}.

  Since $\oM_{g,A}(V,\beta)$ and $\oM_{g,A\simp}(V,\beta)$ are locally
  isomorphic and the second stack is an algebraic Deligne-Mumford
  stack (resp. Artin stack with finite stabilizer), then so is the
  first.

  There are only finitely many equivalence relations, so
  $\oM_{g,A}(V,\beta)$ is of finite type.
  Finally, $\oM_{g,A}(V,\beta)$ is proper by
  Theorem~\ref{thm:properness}.
\end{proof}

The second proof can be applied in the much more general
situation of higher-dimensional stable pairs. We do not discuss here
the many technical complications appearing in the higher-dimensional case.

For every family of stable maps $f:(C,s_i)\to V\times S$, define a
functor
$$\Aut(f): (S\text{-schemes}) \to (\text{Groups})^{\rm op} $$
by setting $\Aut(f)(S')$ to be the automorphism group of
$f':(C',s_i')\to V\times S'$, where $(C',s'_i)=(C,s_i)\times_S S'$.

\begin{theorem}\label{thm:finite-aut}
  $\Aut(f)$ is represented by a finite group scheme over $S$.
\end{theorem}
\begin{proof}
  Locally on $S$ the family $C_U\to U\subset S$ can be embedded into a
  projective space $\bP_U$. By the theory of Hilbert schemes, this implies that
  $\Aut(f)$ is represented by a quasiprojective scheme $G$ over
  $S$. Theorem~\ref{thm:properness} implies that $G\to S$ is proper.

  Finally, it is easy to see that the automorphism group of a stable
  map (over a field) is finite (the weights only help).
  Hence, $G\to S$ is finite.
\end{proof}

\begin{proof}[Second proof of Theorem~\ref{thm:stack-exists}]
  We first show that the maps $f$, together with the extra data of an
  embedding into a fixed projective space, are parameterized by a
  locally closed subscheme in a Hilbert scheme. This is done by a
  standard argument, similar to \cite[Prop. 5.1]{Mumford_GIT}. We then
  take the quotient by $\PGL$ to get rid of the extra data.

  By Remark~\ref{rem:ratnl-weights}, we can assume that each $a_i$
  is
  rational. Let $k$ be a positive integer such that $ka_i\in\bZ$, and
  let $L=\omega_C^k\left(\sum ka_i s_i\right)\otimes f^*\cO_V(1)$. By
  \cite[Prop.3.3]{Hassett_WeightedStableCurves}, for any stable map,
  the sheaf $L^N$ is very ample for $N\ge3$ and does not have higher
  cohomology.

  Let $d=\deg L$ and $M+1=h^0(C,L^N)$. Every choice of a basis in
  $H^0(C,L^N)$ gives an embedding into a fixed projective space
  $\bP^M$ as a closed subscheme with Hilbert polynomial $h(x)=dx +
  1-g$. Further, $Y\times \bP^M$ has a natural ample line bundle
  $p_1^*\cO_Y(1) \otimes p_2^*\cO_{\bP^M}(1)$ of degree $d'=d+
  \deg\beta$. Let $h'(x) = d'x+1-g$ be the Hilbert polynomial
  corresponding to this embedding.

  Now, let $H$ be the Hilbert scheme parameterizing all closed
  subschemes $Z\subset Y\times \bP^M$. The following conditions,
  applied in a sequence, are open: $p_2:Z\to \bP^M$ is a closed
  embedding, $Z$ is reduced, $Z$ has at most nodes, and the
  homomorphism $\Pic V\to \bZ$ induced by $Z\to V$ is $\beta$ (resp.
  $\deg\beta$ is fixed). Let $H_1\subset H$ be the open subscheme
  parameterizing $Z$ with these conditions.

  Next, we add points. First, we choose a closed subscheme $H_2\subset
  H_1 \times (\bP^M)^n$ corresponding to the collections
  $(C,s_1,\dotsc, s_n)$ with $s_i\in C$. Secondly, $H_3\subset
  H_2$ is the open subscheme corresponding to collections such that
  the $n$-pointed curve satisfies the condition on singularities in
  Definition~\ref{defn:stable-map}, which is open in families.

  On the universal family $C_3\to H_3$, we have two invertible sheaves:
  $p_2^*\cO_{\bP^M}(1)$ and $\omega_{C_3/H_3}^{\otimes k}(\sum ka_i
  s_i)$. These give two sections of the relative Picard scheme of
  $C_3/H_3$, which is represented by an algebraic space (see e.g.
  \cite[Thm. 8.3.1]{BoschLutkebohmertRaynaud90}). Let $H_4\subset H_3$
  be the locally closed subscheme where the two sections coincide.

  The morphism $H_4\to \oM_{g,A}(V,\beta)$ is, locally in Zariski
  topology, a $\PGL_{M+1}$-torsor. Indeed, the only difference
  between the two stacks is the embedding of the family into
  $\bP^M$. The stabilizer of the $\PGL$-action is proper by
  Theorem~\ref{thm:finite-aut}. Hence, $\oM_{g,A}(V,\beta)$ is an
  algebraic Artin stack with finite stabilizer.

  The order of the automorphism group of $f$ can be bounded
  universally in terms of $g,n$ and $\deg\beta$. For a stable map defined
  over a field of positive characteristic larger than this bound or a
  field of characteristic zero, the automorphism group scheme is
  reduced. Hence, for large enough divisible $N$, the stack
  $\oM_{g,A}(V,\beta)\times\bZ[1/N]$ is Deligne-Mumford.
\end{proof}

\section{Chambers, walls, and simplicial complexes}
\label{sec:chambers-and-walls}

\begin{definition}
  Let us fix $g,n,V$ and $\beta$.
  We will call the weight data $A\in [0,1]^n$
  \defn{admissible} if the stack $\oM_{g,A}(V,\beta)$ is nonempty.

  Weight data of the same length have a natural partial order:
  $A=(a_i)\ge B=(b_i)$ if $a_i\ge b_i$ for all $i$.
\end{definition}

Following \cite{Hassett_WeightedStableCurves}, we define the set
$\cD_{g,n,\beta}$, where the admissible weights can theoretically
live, and two decompositions of it into finitely many chambers.

\begin{definition}
  If $(g,\beta)\ne (0,0)$ or $(1,0)$, we set $\cD_{g,n,\beta}$ to be the
  cube $[0,1]^n$.  $\cD_{1,n,0}= [0,1]^n$ minus the point $(0^n)$, and
  $\cD_{0,n,0}$ is the subset of $[0,1]^n$ defined by the inequality
  $\sum_{i=1}^n a_i >2$.

We define two chamber decompositions of $\cD_{g,n,\beta}$
obtained by cutting it by finitely many hyperplanes. For the
\defn{coarse
  decomposition}, the hyperplanes are
\[
\cW_c=\left\{ \sum_{i\in I} a_i=1: I\subset \{1,\ldots,n
\},3\le|I|\leq n \right\}. \] For the \defn{fine decomposition},
they are
\[
\cW_f=\left\{\sum_{i\in I} a_i=1: I\subset \{1,\ldots,n
\},2\le|I|\leq n \right\}. \]
\end{definition}

\begin{remark}
We corrected an obvious typo in \cite{Hassett_WeightedStableCurves}
which states $|I|\le n-2$ in place of $|I|\le n$ (which is correct
for $(g,\beta)=(0,0)$).
\end{remark}

There are many kinds of chambers one can define: open, closed, locally
closed. For us, the \emph{locally closed} chambers will be the most useful:
\begin{definition}
  A \defn{chamber} is a nonempty locally closed subset of
  $\cD_{g,n,\beta}$ obtained by choosing for each $I\subset
  \{1,\dotsc,n\}$ either the inequality $\sum_{i\in I} a_i>~1$ or the
  inequality $\sum_{i\in I} a_i\le1$.
\end{definition}

\begin{proposition}[cf. \cite{Hassett_WeightedStableCurves},
Prop.5.1]\label{proposition:chamberconst} For two weight data $A$ and $B$, we have the following:
  \begin{enumerate}
  \item If $A$ and $B$ belong to the same chamber in the coarse
    decomposition and have the same zero weights,
    then $\oM_{g,A}(V,\beta)$ and $\oM_{g,B}(V,\beta)$ naturally coincide.
  \item If $A$ and $B$ belong to the same chamber in the fine
    decomposition and have the same zero weights
    then the universal families over $\oM_{g,A}(V,\beta)$ and
    $\oM_{g,B}(V,\beta)$ naturally coincide.
  \end{enumerate}
\end{proposition}
\begin{proof}\quad
(1) We have to compare the two conditions of
Definition~\ref{defn:stable-map} for $A$ and $B$. The condition on
singularities $\sum_{s_i=P}a_i\le 1$ is obviously the same for $A$ and
$B$. The numerical condition says that for any irreducible component
$E$ of $C$ collapsed by $f$ to a point, one should have
\begin{displaymath}
  \deg \omega_C\left(\sum a_i s_i\right)\Big{|}_E =
  2p_a(E)-2 + |E\cap (C-E)| + \sum_{s_i\in E} a_i >0.
\end{displaymath}
If $(g,\beta)\ne(1,0)$ or $(0,0)$, and $E$ is not a $\bP^1$ with
$|E\cap(C-E)|=2$
this condition is either vacuous or says
\begin{displaymath}
   -1 + \sum_{s_i\in E} a_i > 0
\end{displaymath}
Hence, it is the same for $A$ and $B$.  If $(g,\beta)=(1,0)$ or
$(0,0)$, the additional cases say simply that $A$ and $B$ should be in
$\cD_{g,n,\beta}$. If $E$ is a $\bP^1$ with $|E\cap(C-E)|=2$ then
the condiion says $\sum_{s_i\in E} a_i>0$. Since $A$ and $B$ have the
same zero weights, this condition is equivalent for $A$ and $B$.

  (2) By Lemma~\ref{lem:weights0and1}, the universal family
  $\cC_{g,A}(V,\beta)$ is the moduli space for the weight data $A\cup
  0$; now apply (1).
\end{proof}

\begin{lemma}\label{lem:0-and-eps-admissible}
  Let $(a_i)$ be a positive weight data and suppose that for all
  $0<\epsilon\ll1$, the weight data $(a_i,\epsilon^m)$ is in the
  interior of the same chamber as $(a_i,0^m)$, and both are in
  $\cD_{g,n,\beta}$. Then $(a_i,\epsilon^m)$ is admissible iff
  $(a_i,0^m)$ is admissible.
\end{lemma}
\begin{proof}
  Let $f:(C,s_i)\to Y$ be a stable map w.r.t $(a_i,\epsilon^m)$ for
  \begin{displaymath}
    0< m\epsilon < \min_I \big(-1+ \sum_{i\in I} a_i  > 0 \big).
  \end{displaymath}
  Then the numerical condition of stability fails for the weight data
  $(a_i,0^m)$ on the $f$-exceptional curves $E=\bP^1$ which intersect
  the rest of $C$ at two points and which do not contain any $s_i$
  with $a_i>0$.

  Let $C'$ be the curve obtained by contracting each such $E$ to a
  point.
  Suppose the whole $C$ is contracted this way, i.e. \emph{every}
  component of $C$ was of this form. But then $g=1$, $\beta=0$ and
  $(a_i)=\emptyset$, so the weight data $(a_i,0^m) \not\in
  \cD_{g,n,\beta}$.

  Hence, $C'$ is indeed a curve. The induced map $f':(C',s_i)\to Y$ is
  obviously stable for $(a_i,0^m)$, so the weight data $(a_i,0^m)$ is
  admissible.

  The opposite direction is easier: just put the $m$ points with
  weight $\epsilon$ at some unused nonsingular points of $C'$.
\end{proof}

\begin{corollary}
  For each $(g,n,V,\beta)$ the set of admissible weight data is a
  union of several chambers.
\end{corollary}

We now introduce a convenient way to label the chambers.

\begin{definition}
  Let us identify every subset $I\subset\{1,\dotsc,n\}$ with a simplex
 $\sigma$ with vertices in $I$; we have $\dim\sigma=|I|-1$.
  The \defn{simplicial complex $\Delta_{A}$ associated to the weight data}
  $A=(a_1,\dots, a_n)$ consists of simplices $\sigma(I)$ such that
  $\sum_{i\in I} a_i\le1$. We will often identify $\sigma$ with
  $\sigma(I)$ and denote the simplex as simply $\sigma$.

  If $\sigma(I')$ is a face of $\sigma(I)$, i.e. $I'\subset I$, then
  $\sum_{i\in I}a_i\le1$ implies $\sum_{i\in I'}a_i\le1$. Therefore,
  $\Delta_{A}$ is indeed a complex in the usual sense.
\end{definition}

\begin{remark}
  It is a strong condition for a complex $\Delta$ to be associated
  with a chamber. For example, the complex $\{12,34,1,2,3,4\}$ is not
  associated to any chamber: the system of inequalities $a_1+a_2\le1$,
  $a_3+a_4\le1$, $a_1+a_3>1$, $a_2+a_4>1$ has no solution.
\end{remark}

\begin{corollary} \label{corollary:univcurvedesc}\quad
  With the zero weights fixed:
  \begin{enumerate}
  \item The moduli space $\oM_{g,\Delta}(V,\beta)$ and universal family
$\cC_{g,\Delta}(V,\beta)$ are well-defined  using any $A$ such that
$\Delta=\Delta_{A}$.
\item For the universal family, one has $C_{g,\Delta}(V,\beta) =
  \oM_{g,\Cone\Delta}(V,\beta)$, where $\Cone\Delta$ is the simplicial
  complex on $n+1$ vertices consisting of $\sigma$ and
  $\sigma\cup\{n+1\}$ for all $\sigma\in\Delta$.
  \end{enumerate}
\end{corollary}

In the following example, we fix notation for several complexes of
which we will make use throughout the remainder of this paper.

\begin{example}\label{example:compexamples} \quad
\begin{enumerate}
\item For the weight data $(1^n),$ $\Delta$ is the disjoint union of $n$
  vertices.
\item For the weight data $(\epsilon^n),$ $\Delta$ contains all
  simplices; i.e. the support $|\Delta|$ is homeomorphic to an $(n-1)$-ball.
\item \label{example:symdef} For the weight data $(1/(r+1),...,1/(r+1))$ with $n$ vertices, the complex is
given by the $r$-skeleton of an $n$-simplex (i.e. every subset of
vertices with $r+1$ or fewer vertices is a face of the complex), and
we will denote this complex $\Delta_{n,r}$. Included in this are the
previous two examples for $r=0$ and $r=n-1$, respectively. We will
reference this for any combination of $g,V,\beta$ for which the
corresponding space is nonempty. See Propositions
\ref{theorem:symdilaton} and \ref{proposition:symstring}.
\item \label{example:permutodef} For the weight data $(1^2,\epsilon^{n-2})$ with $\epsilon<1/(n-2)$, the complex $\Delta_{\cL_{n-2}}$ has the
first two vertices isolated, and the remaining vertices form a
complete simplex. Whenever we reference this complex, we will assume
that $g=0,V=\{pt\},\beta=0$. See Lemma \ref{lemma:permuto} and
Example \ref{example:genfunc}.\ref{example:permuto}.
\item \label{example:projdef} For the
weight data $(1,a_{1},\dots,a_{r})$ such that $\sum a_i>1$, but for
any proper subset one has $\sum a_i\le1$, the complex
$\Delta_{\bP^{r-2}}$ is a single isolated vertex plus the boundary
of an $(r-1)$-simplex. Whenever we reference this complex we will
assume that $g=0,V=\{pt\},\beta=0$. See Lemmas \ref{lem:Pr-2},
\ref{lemma:projpsi} and Example
\ref{example:genfunc}.\ref{example:projgen}.

\item Corollary \ref{corollary:univcurvedesc} gives that the complex
  associated to the universal family of $\oM_{g,A}(V,\beta)$ is the
  cone over $\Delta_A$ with vertex given by the additional section. We
  denote this as $\Cone(\Delta_{A})$
\label{example:univcomplex}

\item For any weight data $A$ and $\sigma$ a collection of vertices in $\Delta_A$, we define
$\Delta_{A_{\sigma}}$ to be the complex attained by replacing the
vertices of $\sigma$ with a single vertex, which we shall label
$\sigma$, and assigning the weight equal to the sum of the weights
in $\sigma$. We note that $\gamma\subset\Delta_A$ corresponds to a
face of $\Delta_{A_{\sigma}}$ iff $\gamma$ is a face of $\Delta_{A}$
and either:
\begin{enumerate}
\item $\gamma\cap\sigma=\varnothing$.
\item $\sigma\subset\gamma$.
In this case, $\gamma$ will correspond to the face of
$\Delta_{A_\sigma}$ containing the vertices in
$\gamma\setminus\sigma$ as well as the vertex which we label as
$\sigma$.
\end{enumerate}
This complex will appear naturally in the wall crossing
formula~\ref{theorem:wallcross}.\label{example:sigmacompdef}
\end{enumerate}
\end{example}

\section{Reduction and forgetful morphisms}
\label{sec:reduction-forgetful}

\begin{theorem}[Reduction morphism]
\label{thm:contraction-morphism}
Fix $g,V,\beta$ and let $A, B\in \cD_{g,n,\beta}$ be two admissible
weight data such that $A\ge B$.  Then there exists a natural reduction
morphism
\[
\rho_{B,A}:{\oM}_{g,A}(V,\beta) \ra {\oM}_{g,B}(V,\beta).
\]
Given a stable map $(C,s_1,\ldots,s_n,f)$ for the weight data $A$,
its image \linebreak $\rho_{B,A}(C,s_1,\ldots,s_n,f)$ is obtained by
successively collapsing components of $C$ along which
$K_{\pi}+a_{i_1}s_{i_1}+\ldots+a_{i_r}s_{i_r}$ fails to be
$f$-ample.
\end{theorem}
\begin{proof}
  The proof is the same as in
  \cite[Thm.4.3]{Hassett_WeightedStableCurves}, which is the
  corresponding statement in the absolute case, with the sheaf
  $\omega^k(\sum ka_i s_i)$ replaced by $\omega^k(\sum ka_i
  s_i)\otimes f^*\cO_V(3)$.
\end{proof}

\begin{theorem}[Forgetful morphism]\label{thm:forgetful-morphism}
  Let $A=(a_i)\in \cD_{g,n,\beta}$ and $A'=(a'_i)\in \cD_{g,m,\beta}$
  be two admissible weight data such that $A'$ is obtained from $A$ by
  dropping the last $n-m$ weights. Then there exists a natural
  morphism $\phi_{A',A}:\oM_{g,A}(V,\beta) \to
  \oM_{g,A'}(V,\beta)$. Moreover, $\phi_{A',A}$ is flat and Gorenstein.
\end{theorem}
\begin{proof}
  By the previous theorem, there exists a natural morphism
  $\oM_{g,A}(V,\beta)$ $\to\oM_{g,A'\cup 0^{n-m}}(V,\beta)$, and
  $\oM_{g,A'\cup 0^{n-m}}(V,\beta) \to \oM_{g,A'}(V,\beta)$ is the
  $n$-th fiber power of the universal curve.

  By Corollary \ref{cor:local-structure},
  locally both
  $\oM_{g,A}(V,\beta)$ and $\oM_{g,A'}(V,\beta)$ are isomorphic to the
  stacks for the weight data $A^R_{\rm sim}$ and $A'^{R'}_{\rm
    sim}$. Here, the weight datas $A$ and $A'$ consist of $1$'s and
  $0$'s only, and the index set of $A'$ is naturally a subset of the
  index set $A$.
  Thus, locally, the moduli spaces are isomorphic to
  $\oM_{g,(1^a,0^b)}(V,\beta)$, $\oM_{g,(1^{a'},0^{b'})}(V,\beta)$
  with $a'\le a$,   $b'\le b$.

  In these local charts, $\phi_{A',A}$ reduces to the natural
  forgetful morphism $\oM_{g,(1^a,0^b)}(V,\beta)\to
  \oM_{g,(1^{a'},0^{b'})}(V,\beta)$,
  which is a composition of morphisms forgetting a single 1 or 0. Each
  of these elementary morphisms is a universal family over the
  corresponding moduli space. Hence this morphism is flat. Since the
  curves are nodal, it is also Gorenstein.
\end{proof}


\begin{lemma}[Decreasing from $\epsilon$ to $0$]
  \label{lem:reduction-to-zero}
  Let $(a_i)$ be a positive weight data and suppose that for all
  $0<\epsilon\ll1$, the weight data $A_{\epsilon}=(a_i,\epsilon^m)$ is
  in the interior of the same chamber as $A_0=(a_i,0^m)$, and the latter
  is admissible (cf. Lemma~\ref{lem:0-and-eps-admissible}).

  For each $[f:(C,s_i)\to Y]\in \oM_{g,A_0}\vb$, suppose some $k_1$
  zero-weight points coincide with one node, and $k_2$ with another
  node, etc.

  Then the fiber of the morphism $\rho: \oM_{g,A_{\epsilon}}\vb \to
  \oM_{g,A_0}\vb$ over $[f]$ is $\prod_{\alpha}\oM_{0,
    (1^2,\epsilon^{k_{\alpha}})}$. In particular, it is a single point
  if and only if all $k_{\alpha}=1$, and in this case $\rho$ is an
  isomorphism over a neighborhood of $[f]$.
\end{lemma}
\begin{proof}
  Indeed, as in the proof of Lemma~\ref{lem:0-and-eps-admissible}, the curves
  that may get identified by $\rho$ are obtained by contracting
  $f$-exceptional $E=\bP^1$ with two points of weight 1 and
  $k_{\alpha}$ points of weight $\epsilon$. Finally,
  $\oM_{0,(1^2,\epsilon)}$ is a point.
\end{proof}

\begin{lemma}[Universal family]
  \label{lem:universal-family}
  Let $A=(a_i)$ be an admissible positive weight data. Then for all
  $0<\epsilon\ll \delta\ll 1$, the universal family over
  $\oM_{g,A}\vb$ is the moduli space for the weight data
  $(a_i-\delta, \epsilon)$.
\end{lemma}
\begin{proof}
By Lemma \ref{lem:weights0and1}, the moduli
  space $\oM_{(b_i,0)}(V,\beta)$ is naturally isomorphic to the
  universal family over the moduli space $\oM_{(b_i)}(V,\beta)$.

  Now, by Remark~\ref{rem:ratnl-weights}, the moduli spaces for
  $(a_i)$ and $(a_i-\delta)$ are the same for $0<\delta\ll1$. By choosing
  $\delta$ sufficiently small, we can guarantee that the point
  $(a_i-\delta,\epsilon)$ is in the interior of the same chamber as
  $(a_i-\delta,0)$, so that the previous lemma applies.

  Since in this situation there is a unique $k_1$, equal to 1, by the
  previous lemma the morphism $\oM_{g,(a_i-\delta,\epsilon)}\vb \to
  \oM_{g,(a_i-\delta,0)}\vb$ is an isomorphism.
\end{proof}

\section{Crossing a single wall = adding a simplex}
\label{sec:crossing-one-wall}

\begin{definition}
  If $\Ap\ge A$ and $\Delta_A$ is obtained from $\Delta_{\Ap}$ by
  adding a single simplex $\sigma(I)$, i.e. by changing the sign in
  the single inequality from $\sum_{i\in I}a_i>1$ to $\le1$, then the
  change from $\Ap$ to $A$ will be called a {\bf simple wall
    crossing}.
\end{definition}

\begin{lemma}
  For any two \emph{positive} weight data $A\ge B$ there exist weight data
  $A'\ge B'$ such that $\Delta_A=\Delta_{A'}$, $\Delta_B=\Delta_{B'}$,
  and the the straight line from $A'$ to $B'$ goes through a sequence
  of simple wall crossings.
\end{lemma}
\begin{proof}
  Indeed, we can change $(a_i)$ by $(a_i-\epsilon_i)$ for small generic
  $\epsilon_i>0$, and similarly for $b_i$. The condition that the line
  passes through an intersection of the hyperplanes
  $\sum_{i\in I}a_i=1$ is a union of hyperplanes, and so for the generic
  $\epsilon_i$ this does not happen.
\end{proof}

\begin{definition}\label{defn:IJ}
Let $J$ be the complement of $I$.  Consider
  two weight data $\Ap=(a_i^+,a_j)$ and $A=(a_i,a_j)$, $i\in I$, $j\in
  J$ such that
\begin{enumerate}
\item $a_i^+\ge a_i>0$ for $i\in I$,
\item $\sum_{i\in I} a_i =1$ and $\sum_{i\in I} a^+_i >1$,
\item $\Delta_{A}$ is obtained from $\Delta_{\Ap}$ by adding a single
  simplex $\sigma(I)$.
\end{enumerate}
In this situation, we will introduce new weight data
 \begin{displaymath}
   A_I^+ = (1,a_i^+) \quad\text{and}\quad A_J=(1,a_j)
 \end{displaymath}
of lengths $|I|+1$ and $|J|+1$ respectively.
We denote $r=|I|$ for simplicity.
\end{definition}

\begin{theorem}\label{thm:complete-intersection}
  There exists a natural closed embedding $\iota:\oM_{g,A_J}(V,\beta)
    \to \oM_{g,A}(V,\beta)$ defined by associating to a family
  $$(C,\sbul,s_{r+1},\dotsc,s_n)\to S, \quad f:C\to V\times S$$
  the family
  $$(C,\sbul,\dotsc,\sbul, s_{r+1},\dotsc,s_n)\to S, \quad \text{ with the same }f$$
  in which the section $\sbul$ is repeated $r$ times.

  The image is an intersection of $r-1$ globally defined Cartier
    divisors which is locally a complete intersection; in particular
    it has codimension $r-1$.
\end{theorem}
\begin{proof}
  Since $\sum_{i\in I} a_i =1$, it is clear that
  the first family in the statement of the theorem is stable with
  respect to the weight data $A_J$ if and only if the second one is
  stable with respect to the weight data $A$. Hence, $\iota$ is well
  defined.

  As a set, the image is the intersection of $r-1$ globally defined
  closed sets $\{s_1=s_i\}$, $i=2,\dots,r$. Let us check that these
  sets are indeed given by Cartier divisors, and that locally it is a
  complete intersection.

  To see this, we can work with a local chart. Thus,
  we can simplify the weight $A$ partially, by replacing
  the part $(a_i)$, $i\in I$, by $(1,0,\dotsc,0)$.  According to
  Corollary~\ref{cor:local-structure}, this will not change the moduli
  stack in a neighborhood of the image of $\iota$.

  By Lemma~\ref{lem:weights0and1}, the moduli space for the partially
  simplified weight is an $(r-1)$st fiber power of the universal curve
  over $\oM_{g,A_J}\vb$, and $\iota$ is a section of this projection.
  Since the morphism $C\to \oM_{g,A_J}\vb$ of the universal family is
  smooth of relative
  dimension one along the section $\sbul$, the morphism from the fiber
  power is smooth of relative dimension $r-1$ along $(\sbul,\dotsc,\sbul)$.
  Hence, its section is a complete intersection of $r-1$
  Cartier divisors $F_i$ on which $s_{i}=\sbul$.
\end{proof}

We now describe the reduction morphism $\rho_{A,\Ap}$. We start
set-theoretically.

\begin{lemma}\label{lem:Pr-2}
  Let $a^+_1,\dots,a^+_r$ be nonnegative numbers such that $\sum_{i=1}^r
  a^+_i>1$, but for any proper subset one has $\sum a^+_i\le1$. Let
  $A_I^+=(1,a^+_i)$ be the weight data of length $r+1$. Then $\oM_{0,A_I^+}
  =\bP^{r-2}$, the universal family is the blowup of $\bP^{r-1}$ at a
  point, 
  and every geometric fiber $\cC_{\bar s}$ in the universal family is
  isomorphic to $\bP^1$.
\end{lemma}
\begin{proof}
  If $C_{\bar s}$ is not $\bP^1$, then it is a tree of $\bP^1$'s with
  at least two endpoints. For the corresponding irreducible components
  $E_j$ of $C$, one must have $\sum_{s_i\in E_j}a^+_i>1$, which is not
  possible by the conditions.

  Let us denote the section corresponding to weight 1 by $\infty$, and
  the section corresponding to $a^+_1$ by $0$. Then $\oM_{0,A_I^+}$ is the
  moduli of $r$ points on $\bA^1$, not all of which are equal to
  $0$. Taking into account that $\Aut(\bP^1,0,\infty)=\bG_m$, we get
  $\oM_{0,A'}= (\bA^{r-1}\setminus 0)/\bG_m = \bP^{r-2}$.

  Consider $\bP^{r-1}$ with $r$ coordinate hyperplanes and a point
  $P=(1,\dotsc, 1)$. Let $\cC=\Bl_P \bP^{r-1}$ be the blowup, and
  $\pi:\cC\to \bP(T_p) = \bP^{r-1}$ be the natural fibration with $r+1$
  sections: $r$ for the hyperplanes and one more for the exceptional
  divisor.
  If we assign to the hyperplanes weights $a^+_i$ and to the exceptional
  divisor weight $1$, then we get the universal family of stable curves
  for the weight data $A_I^+$ over $\bP^{r-2}$.

  Indeed, consider the family of lines through the point
  $(1,\dotsc,1)$. We can write the equation of such a line
  parametrically as follows: $t_1 x_i = x_1 (t_1- c_i t_0)$,
  $(t_1:t_0)$ being a point of $\bP^1$. In this parameterization,
  $(1,\dotsc,1)$ corresponds to $(t_0:t_1)=(0:1)$, the
  intersection with the hyperplane $x_1=0$ to the point $(1:0)$, and
  the coefficients $c_i$ are defined up to a common invertible
  multiple. The intersections with the hyperplanes $x_i=0$,
  $i=2,\dotsc,r$ are the points $(1:c_i)$.  Hence, this is the
  universal family over $(\bA^{r-1} \setminus 0)/\bG_m = \bP^{r-2}$.
\end{proof}

\begin{lemma}\label{lem:DIJ}
  Let $p$ be a geometric point of $\oM_{g,A}(V,\beta)$. Then
  \begin{enumerate}
  \item if $p\not\in \iota\big(\oM_{g,A_J}(V,\beta)\big)$,
    then $\rho_{A,\Ap}\inv(p)$ is one point corresponding to the same
    map $f$.
  \item if $p\in \iota\big(\oM_{g,A_J}(V,\beta)\big)$
    then $\rho_{A,\Ap}\inv(p)$ is $\bP^{r-2}$ and consists of
maps $f:C\to V$ such that $C=C_I \cup
C_J$ is a union  of two curves $C_I,C_J$ glued at the point $\sbul$
for the extra weight 1, such that $C_I=\bP^1$ is collapsed to
a point of by $f$, and $f|_{C_J}$ is a stable map for the weight data
$A_J$.
  \end{enumerate}
  Set-theorically, on the exceptional sets the morphism is the
  projection to the second summand
  $$ \bP^{r-2} \times \oM_{g,A_J}(V,\beta) \to \oM_{g,A_J}(V,\beta)
  \to \iota   \oM_{g,A_J}(V,\beta) $$
\end{lemma}
\begin{proof}
  It is easy to see that the only contraction that can occur is a
  $\bP^1$ which intersects the rest of the curve $C$ at one point,
  contains precisely the points $s_i$ for $i\in I$, and $C_I$ is
  collapsed by $f$. Hence, $C=C_I \cup C_J$ as claimed. The set of
  such curves is the moduli space of $(C_I, s_i, \sbul)$, which is
  $\bP^{r-2}$ by the previous lemma.
\end{proof}

\begin{theorem}\label{thm:locally-smooth}
  There exists an open neighborhood $U$ of
  $\iota\big(\oM_{g,A_J}(V,\beta)\big)$ with a morphism $\phi':U\to
  \oM_{g,A_J}(V,\beta)$ such that $\phi'$ and $\phi'\circ\rho_{A,\Ap}:
  \rho_{A,\Ap}\inv U\to \oM_{g,A_J}(V,\beta)$ are smooth.
\end{theorem}
\begin{proof}
  The smooth morphism $\phi'$ was constructed in the proof of
  Theorem~\ref{thm:complete-intersection}. The morphism
  $\phi'\circ\rho_{A,\Ap}$ is a restriction of a forgetful morphism to
  an open set. By Theorem~\ref{thm:forgetful-morphism} we know that it
  is flat, and that it can be decomposed into several morphisms from
  the universal family to the corresponding moduli stack. Such a
  ``universal'' morphism is
  smooth except at the node. When one encounters a node, the new
  moduli space (the universal family over an old moduli) corresponds
  to some curves with 2 new nodes. This does not happen by
  Lemma~\ref{lem:DIJ}. Hence, $\phi'\circ\rho_{A,\Ap}$ is smooth.
\end{proof}

We now describe the reduction morphism $\rho_{A,\Ap}$ as an explicit
blowup along a complete intersection. This is true even though the
moduli spaces involved could be very complicated: nonreduced and not
equidimensional.  Philosophically, such a nice description is
possible because, by the above theorem, the ``hard part'' of the
moduli -- varying the curve and the map -- and the ``easy part'' --
varying points -- locally split. The change from $\Ap$ to $A$ occurs
only in the ``easy'' direction.

\begin{theorem}[cf. \cite{Hassett_WeightedStableCurves}, Prop.4.5]
\label{thm:blowup}
For a simple wall crossing,\linebreak
$\oM_{g,\Ap}(V,\beta)$ is the blowup of $\oM_{g,A}(V,\beta)$ along
$\oM_{g,A_J}(V,\beta)$, which is a complete intersection of
codimension $r-1$, where $r:=|I|$.
The exceptional Cartier divisor is
\[
D_{I,J}:=\oM_{0,\Ap_I} \times\oM_{g,A_J}(V,\beta), \quad
\text{moreover} \quad \oM_{0,\Ap_I} = \bP^{r-2},
\]
which as a set was described in Lemma~\ref{lem:DIJ}.
\end{theorem}
\begin{proof}
  Let $\cI$ be the ideal of $\iota\big(\oM_{g,A_J}(V,\beta)\big)$ in
  $\oM_{g,A}(V,\beta)$. We first prove that $\rho=\rho_{A,\Ap}$ factors
  through the blowup:
  \begin{displaymath}
    \rho':   \oM_{g,\Ap}(V,\beta)\to
  \operatorname{Bl}_{\iota(\oM_{g,A_J}(V,\beta))} \oM_{g,A}(V,\beta)
  \end{displaymath}

  For this, we check the universal property of
  the blowup: the preimage $\rho\inv \cI$ of the ideal sheaf $\cI$ of
  $\oM_{0,A_J}\vb$ is invertible.

  From the description given in the proof of
  Theorem~\ref{thm:complete-intersection}, the ideal $I$ is locally generated
  by regular functions $f_2, \dots, f_r$ such that the zero set of
  $f_i$ is the locus in $\oM_{0,A}\vb$ where the points $s_1$ and $s_i$
  coincide. Let $g$ be a local equation of the exceptional divisor
  $E$. Then we have:
  \begin{math}
    \rho^*(f_i) = gf'_i,
  \end{math}
  and the zero set of $f'_i$ is contained in the locus in $\oM_{0,\Ap}\vb$
  where the points $s_1$ and $s_i$ coincide. Hence, the zero set of
  the  ideal $(f'_i)$  is contained in the locus in $\oM_{0,\Ap}\vb$ where
  all $r$ points $s_i$ coincide. This locus is empty. Therefore,
  \begin{displaymath}
    \rho\inv \cI = (gf'_i) = (g)
  \end{displaymath}
  is principal, and so $\rho$ factors through the blowup.

  On the other hand, by Theorem~\ref{thm:complete-intersection}, we
  know that $\oM_{g,A_J}(V,\beta)\subset \oM_{g,A}(V,\beta)$ is a
  complete intersection of codimension $r-1$, so its normal bundle is
  a direct sum of $r-1$ copies of the same line bundle. Therefore, the
  exceptional set of the blowup is also $\bP^{r-2}\times\oM_{0,A_J}\vb$.

  Hence, the morphism $\rho'$ to the blowup is a bijection on
  geometric points. Since both the source and the target of $\rho'$
  are smooth over $\oM_{g,A_J}\vb$ by
  Theorem~\ref{thm:locally-smooth}, it follows that $\rho'$ is an
  isomorphism.
\end{proof}

\begin{corollary}\label{cor:blowupcor}
  For any two positive weight data $A\ge B$, the stable reduction
  morphism $\rho_{B,A}:\oM_{g,A}\vb\to \oM_{g,B}\vb$ is a composition of
  several blowups; it contracts the divisors $\Ds$ associated to the
  faces $\sigma$
  of dimension strictly greater than one which are in the complex
  $\Delta_B$ but not in $\Delta_A$.

  The morphism $\rho_{B,A}$ is an isomorphism iff $\Delta_A$ and
  $\Delta_B$ differ only in several edges, i.e. simplices of dimension
  1.
\end{corollary}

\section{Psi classes}

Let $\pi_{A}:\cC_{g,A}(V,\beta)\to\oM_{g,A}(V,\beta)$ be the map of
the universal family with sections $s_{i,A}$ and relative dualizing
bundle $\omega_{\pi_{A}}$, and let $N_{s_i}$ be the
normal bundle of $s_i$ in the universal family.

\begin{definition}\label{defn:psi-classes}
For a \emph{positive} weight data $A$, we define the \defn{psi
classes} to be
\[
\psi_{i,g,A}:=c_1 \big({s_{i,A}^*(\omega_{\pi_{A}})} \big)
= - c_1(N_{s_i})
\]

Note that if $a_i>0$, then by Definition~\ref{defn:stable-map},
$s_i$ is contained in the locus of $\pi:C\to S$ where $\pi$ is
smooth. Hence, $\psi_{i,g,A}$ are the first Chern classes of
invertible sheaves. If some of the weights $a_i=0$, we must adjust
this definition.

Let $A= (a_i)$ be the positive weight data. Then by
Lemma~\ref{lem:universal-family}, the universal family $\cC$ over
$\oM_{g,A}\vb$ is $\oM_{g,(a_i-\delta,\epsilon)}\vb$, and so for the
weight-$\epsilon$ section there is a well-defined psi class using
the above definition, and let us denote it simply $\psi$.

For the weight data $A\cup 0^m$ we define the psi class for the
section $s_j$ with $a_j=0$ to be the pullback of $\psi$ under the
$j$-th projection
\begin{displaymath}
  \oM_{g,A\cup 0^m}\vb  = \cC^m \to \cC: \qquad
  \psi_j: = p_j^* (\psi).
\end{displaymath}
Here, $\cC^m$ is the $m$-th fibered power over $\oM_{g,A}\vb$.
\end{definition}

We may refer to the psi class of a vertex of $\Delta$ with the
obvious meaning. We also omit any subscripts of the notation
whenever it is unlikely to lead to confusion.

\begin{lemma}\label{lemma:projpsi}
Consider the complex $\Delta_{\bP^{r-2}}$ as defined in Example
\ref{example:compexamples}.\ref{example:projdef}. The psi classes of
the nonisolated vertices are $-h$, for $h$ the hyperplane section of
$\bP^{r-2}$, and the isolated vertex has psi class $h$.
\end{lemma}
\begin{proof}
Consider the map $\pi$ of the universal curve, which is described in
Lemma \ref{lem:Pr-2}. For the preimages of hyperplanes $H_i$, resp.
for the exceptional divisor $E$ of the blowup, we get
\begin{displaymath}
N_{s_i} = \cO_{H_i}(H_i) = \cO_{H_i}(1),
\quad \text{resp.} \quad
N_{s_i} = \cO_{E}(E) = \cO_E(-1).
\end{displaymath}
Hence, $\psi_i = -c_1(N_{s_i})$ are as claimed.
\end{proof}

\subsection{Pull-back via reduction}
Recall that $D_\sigma$ is the divisor parameterizing maps of curves
with sections corresponding to $\sigma$ on a contracted genus zero
component and the remaining sections on the genus $g$ component as
in Corollary \ref{cor:blowupcor}.

\begin{theorem}\label{thm:pullpsi}
For any simple wall
  crossing with positive weight data $A$, we have
$\psi_{i,\Ap}=\rho^*(\psi_{i,\Am})+\Ds$ for $i\in \sigma$ and
$\psi_{j,\Ap}=\rho^*(\psi_{j,\Am})$ for $j\not\in\sigma.$
\end{theorem}

\begin{proof}
  Let $\rho=\rho_{A,\Ap}$ be the reduction morphism between the moduli
  stacks, and $\widetilde{\rho}=\rho_{A\cup 0, \Ap\cup 0}$ be the induced reduction
  morphism between their universal families. For clarity, we denote
  $\Ds=D_{I,J}$. Also, recall the notations introduced in
  Definition~\ref{defn:IJ}.

  By choosing some $0< \epsilon \ll \delta\ll 1$, we can replace
  $A=(a_i)$ by the weight data $A'=(a_i-\delta)$ without changing the
  chamber, so that $A'\cup 0$ and $A'\cup\epsilon$ lie in
  the same chamber as well; we will keep denoting it $A$. Similarly, we can assume that the
  universal family for $\Ap$ is the moduli stack for the weight
  $\Ap\cup \epsilon$.
  Then we see that $\widetilde{\rho}$ is a composition of two simple wall
  crossings: for $(I\cup\bullet,J)$ and for $(I,J\cup\bullet)$.

We get the following commutative diagram:
\begin{equation*}
\xymatrix{
\oM_{g,\Ap\cup\epsilon}(V,\beta) \ar[d]^{\pi_{A^+}}
\ar[rr]^{\widetilde{\rho}} &&
\oM_{g,A'\cup\epsilon}(V,\beta) \ar[d]^{\pi_A} \\
\oM_{g,\Ap}(V,\beta) \ar[rr]^{\rho} &&
\oM_{g,A'}(V,\beta)
}
\end{equation*}


Now we would like to write the standard formulas for the change of the
canonical class under the blowup. This would work if the above moduli
stacks were smooth, but of course they may be highly singular,
not equidimensional, etc.

So, instead we note that all of these spaces come with forgetful
morphisms to the moduli space $\oM_{J}:=\oM_{g,A_J}(V,\beta)$, and by Theorem
\ref{thm:forgetful-morphism} they are all Gorenstein over $\oM_{J}$.
Therefore, the relative dualizing sheaves over $\oM_J$ exist and are
invertible. We have:
  \begin{eqnarray*}
    \omega_{\oM_{\Ap}/\oM_J} &=& \rho^* \omega_{\oM_{A'}/\oM_J}
    + (r-2) D_{I,J} \\
    \omega_{\oM_{\Ap\cup\epsilon}/\oM_J} &=&
    \widetilde{\rho}^* \omega_{\oM_{A'\cup\epsilon}/\oM_J}
    + (r-2) D_{I,J\cup\bullet}
    + (r-1) D_{I\cup\bullet,J}
  \end{eqnarray*}
  Subtracting, and taking into account $\pi_{\Ap}^*(D_{I,J}) =
  D_{I,J\cup\bullet} + D_{I\cup\bullet,J}$, we obtain
  \begin{displaymath}
    \omega_{\pi_{\Ap}} = \widetilde{\rho}^* \omega_{\pi_A} + D_{I\cup\bullet, J}
  \end{displaymath}
  The restriction of $D_{I\cup\bullet, J}$ to a section $s_i,i\in
  I$ is $D_{I,J}$, and to a section $s_j,j\in J$ is 0. This gives the
  stated formula.
\end{proof}

\begin{corollary}\label{corollary:pullpsi}
  For any positive $A\ge B$, let $F(A,B)$ be the set of faces of
  $\Delta_B$ which are not in $\Delta_A$. Then we have
\[
\psi_{i,A}=\rho_{B,A}^*(\psi_{i,B})+\sum_{\stackrel{\sigma\in
F(A,B)}{i\in\sigma}}\Ds
\]
\end{corollary}

\subsection{Pullback via decreasing from $\epsilon$ to $0$}

The next lemma pertains to one of the first examples of weighted
moduli spaces which were studied by A. Losev and Yu. I. Manin in
\cite{LosevManin}. It is, interestingly enough, the toric variety
associated to the permutohedron, the convex hull of the
$S_n$--orbit of $(1,2,\dots ,n)$. In fact, this moduli space can be
interpreted as the moduli space of stable $(n-2)$-pointed curves
$(\mathbb G_m\curvearrowright C,s_2,\dotsc,s_{n-1})$ with
torus action.

\begin{lemma}\label{lemma:permuto}
  Consider the complex $\Delta_{\cL_{n-2}}$ for the weight data
  $(1^2,\epsilon^{n-2})$, as in Example
  \ref{example:compexamples}.\ref{example:permutodef}. Then the psi
  classes of the nonisolated vertices are zero.
\end{lemma}
\begin{proof}
  Decrease $a_2=1$ to $1-(n-3)\epsilon$. This gives
  the reduction morphism $\oM_{\cL_{n-2}} \to
  \oM_{\bP^{n-3}}$. On the latter space, which is isomorphic to
  $\bP^{n-3}$, the psi class is $-h$. Under the
  successive wall crossings, the psi class is changed by adding
  exceptional divisors for blowing up at $n-3$ points, then $\binom
  {n-3} 2$ strict preimages of lines though those points, then
  $\binom{n-3} 3$ strict preimages of 2-planes, etc., and ending with
  a strict preimage of a hyperplane. This adds up to the full preimage
  of a hyperplane, and the result is zero.
\end{proof}

  (See also a second proof at the end of Subsection~\ref{section:pullforget}.)

\begin{theorem}\label{thm:psi-eps-to0}
  As in Lemma~\ref{lem:reduction-to-zero}, let $(a_i)$ be the positive
  weight data and suppose that for all
  $0<\epsilon\ll1$, the weight data $A_{\epsilon}=(a_i,\epsilon^m)$ is
  in the interior of the same chamber as $A_0=(a_i,0^m)$, and the latter
  is admissible. Let $\rho:\oM_{g,A_{\epsilon}}\vb \to \oM_{g,A_0}\vb$
  be the reduction morphism. Then
  \begin{displaymath}
    \psi_{i,A_{\epsilon}} = \rho^*\psi_{i,A_0} \quad\text{for all } i.
  \end{displaymath}
\end{theorem}
\begin{proof}
  By Lemma~\ref{lem:reduction-to-zero} the fibers of $\rho$ are
  products of $\oM_{0, (1^2,\epsilon^k)}$, and by
  Lemma~\ref{lemma:permuto} the restrictions of
  $\psi_{i,A_{\epsilon}}$ to the fibers are all zero. Together with
  the fact that by Theorem~\ref{thm:forgetful-morphism} both spaces
  are Gorenstein over $M_{g,(a_i)}\vb$, this gives the statement.
\end{proof}

\subsection{Pull-back via the forgetful morphism}
\label{section:pullforget} We start by considering the pullback of
psi classes whenever the forgetful map corresponds to that of the
universal curve as given in Lemma \ref{lem:weights0and1} and
Corollary \ref{corollary:univcurvedesc}.

\begin{lemma}\label{lem:drop0}
  Let $A$ be any admissible weight data of length $n$, and let
  \linebreak $\phi:
  \oM_{g,A\cup 0}\vb \to \oM_{g,A}\vb$ be the forgetful morphism. Then
  \begin{displaymath}
    \psi_{i,A\cup 0} = \phi^*\psi_{i,A}
    \quad\text{for all } 1\le i\le n.
  \end{displaymath}
\end{lemma}
\begin{proof}
  Indeed, the universal family over $\oM_{g,A\cup 0}\vb$ is the
  cartesian product
  \begin{displaymath}
    \cC_{g,A}\vb \times_{\oM_{g,A}\vb} \oM_{g,A\cup 0}\vb,
  \end{displaymath}
  and so $\omega_{\pi\cup 0}$ is the pullback of $\omega_{\pi}$.
\end{proof}

We recall that in the unweighted case, we have the well known basic
pullback relationship due to Witten which states that
\[
\psi_{i,n+1}=\phi^*(\psi_{i,n})+D_{i,n+1},
\]
where $D_{i,n+1}$ is the divisor with only the marked points $i,
n+1$ on a genus zero contracted component, and $\phi$ is the
morphism which forgets the $n+1^{st}$ point and stabilizes.

Define $A':=A\setminus\{a_{n+1}\}$ and consider the complexes
$\Delta_{A}$ and $\Cone(\Delta_{A'})$. We identify the vertex of
$\Cone(\Delta_{A'})$ with the last section of $A$. We note that,
just as in the standard unweighted case, the faces
\[
F(A,A'):=\{\sigma\in\Cone(\Delta_{A'}): \sigma\not\in\Delta_A\}
\]
are in bijection with divisors $D_{\sigma}$ which become unstable
after forgetting the last section. Using these observations, we
state and prove the analogue of the basic pullback relationship.

\begin{theorem}[Basic Pullback Relation]
  \label{thm:forback}
  For any admissible weight data $A$ of length $n+1$, if $A':=A\setminus\{a_{n+1}\}$ is
  also admissible, then
\[
\psi_{i,A}= \phi_{A',A}^*(\psi_{i,A'})+\xi_{i}\\
\]
with
\[
\xi_{i}:=\sum_{\stackrel{\sigma\in F(A,A')}{i\in\sigma}} \Ds
\]
\end{theorem}

\begin{proof}
  By Lemma~\ref{lem:drop0} we can assume that $A$ is positive. We
  decrease the weights $a_i$ a little to get into the interior of a
  chamber. Then we decrease the weight $a_{n+1}$ we are about to
  forget to $\epsilon$, and then to $0$. The psi classes will change
  as claimed by Corollary~\ref{corollary:pullpsi} and
  Theorem~\ref{thm:psi-eps-to0}. Then we apply Lemma~\ref{lem:drop0}
  one more time.
\end{proof}

We use this result to give the second proof of
Lemma~\ref{lemma:permuto} now.

\begin{proof}[Second proof of Lemma~\ref{lemma:permuto}]
  By Theorem~\ref{thm:forback}, the psi classes of the nonisolated vertices
  pull back from $\Delta_{\cL_{1}}$. $\oM_{\Delta_{\cL_{1}}}$
  is isomorphic to $\oM_{0,3}$ = a point. Therefore, the
  pullback is zero.\\

\end{proof}

\section{Gravitational Descendants}

A crucial ingredient in the theory of stable maps is the notion of the
virtual fundamental class of $\oM_{g,A}(V,\beta)$. In the unweighted
case, this is treated in \cite{BehrendFantechi} and \cite{Behrend},
among other places. In the weighted context, it is reasonable to expect that
one could define a virtual fundamental class in the same fashion.
Rather than repeat those constructions, we make the following
definition of the virtual fundamental class for the weight $A$.

\begin{definition}[Virtual Fundamental Class]\label{defn:virtual-class}
Let $A$ be an admissible weight data and $\rho$ the reduction
morphism from $(1^n)$ to $A$. Then we define
\[
[\oM_{g,A}(V,\beta)]\virt:=\rho_*[\oM_{g,n}(V,\beta)]\virt
\]
\end{definition}
We do not check the
axioms of the virtual fundamental class, which in the unweighted case
were given in \cite{BehrendManin}; for this, one
may consult \cite{BayerManin}.

The moduli spaces $\oM_{g,A}(V,\beta)$ are equipped with $n$
evaluation morphisms $\nu_{i,A}:\oM_{g,A}(V,\beta)\to V$ defined by
$\nu_{i,A}\Big{(}[C,\{s_i\}, f]\Big{)}=f(s_i).$

\begin{definition}\label{defn:descendants}
We define an
analogue of the usual notion of the gravitational descendants of
Gromov-Witten theory which we denote as
\[
\des{\tau_{k_1}(\gamma_1)\cdots\tau_{k_n}(\gamma_n)}{V,\beta}{g,A}:=
\int \left(\prod\psi_{i,A}^{k_i}\cup\nu^*_{i,A}(\gamma_i)\right)
\bigcap\Big{[}{\oM_{g,A}(V,\beta)}\Big{]}\virt,
\]
where $\gamma_i\in A^*(V,\bQ)$ ($\gamma_i\in H^*(V,\bQ)$ when working
over $\bC$), and each $k_i$ is a nonnegative integer.

As usual, these are defined to be zero unless
\[
\sum_{i=1}^{n}(k_i+\deg\gamma_i)=(1-g)\dim V -K_V\beta + (3g-3+n).
\]
Whenever any $k_i$ is negative, we define this to be zero as well.
\end{definition}
We note that whenever $k_i=0,$ these are simply the Gromov-Witten
invariants of $V.$

We warn the reader to treat the $\tau$'s as noncommuting variables
and to not shift indices without discretion as the symmetry of these
descendants is very often broken. One can describe the commuting
properties of the $\tau$'s in terms of the symmetries of the complex
$\Delta_A$, but we will make no use of this and leave it to the
reader. When the weight is $(1^n),$ we omit the weight and note the
number of marked sections, the genus $g$, $V$ and $\beta$.

The first property of the weighted descendants is this:

\begin{lemma}
  For $g,n,V,\beta$ fixed, each descendant
  $\des{\tau_{k_1}(\gamma_1)\cdots\tau_{k_n}(\gamma_n)}{V,\beta}{g,A}$
  is constant as $A$ varies in a fine chamber.
\end{lemma}
\begin{proof}
  Indeed, when $A$ varies so that the zero weights remain the same,
  the moduli space and the universal family stay constant by
  Proposition~\ref{proposition:chamberconst}. When some coefficients
  decrease from $\epsilon$ to $0$, the intersection stays constant by
  Theorem~\ref{thm:psi-eps-to0}, Definition~\ref{defn:virtual-class}
  and projection formula.
\end{proof}

Some weighted descendants are actually very familiar:

\begin{lemma}\label{lem:kappas}
The descendant for the weight data $(0^n)$ are the
intersections of the Miller-Morita-Mumford kappa classes
\[
\des{\tau_{k_1}\cdots\tau_{k_n}}{}{g,(\epsilon^n)} =
\des{\tau_{k_1}\cdots\tau_{k_n}}{}{g,(0^n)} =
\des{\kappa_{k_1-1}\cdots\kappa_{k_n-1}}{}{}
\]
\end{lemma}
\begin{proof}
The first identity is by
  Theorem~\ref{thm:psi-eps-to0}.
  The second identity is simply by the definition of the kappa classes
  (see e.g. \cite{Witten}) and by our definition of the psi classes on
  for the weight data $(0^n)$.
\end{proof}

\begin{subsection}{The Descendant Invariants of \cite{GKP}}\label{subsection:gkp}

In \cite{GKP}, a modification of the psi classes is defined which we
reinterpret in the weighted context.

Suppose $\beta>0$ or $g>0$.
    For each mark $p_i$, let

\[
    \hat\pi_i :\oM_{g,n}(V,\beta) \to \oM_{g,\{p_i\}}(V,\beta)
\]
    be the morphisms which forget all marks but $p_i$.   The modified
    psi class on $\oM_{g,n}(V,\beta)$ is by definition
\[
   \overline{\psi}_i := \hat\pi_i^* (\psi_i).
\]
Therefore, these modified psi classes are simply the pullbacks of the
psi classes for the weight data $(0^n)$, as we defined them above, via
the reduction morphism $\rho:\oM_{g,(1^n)}(V,\beta) \to
\oM_{g,(0^n)}(V,\beta)$. The modified gravitational descendants using
psi classes $\psi_1, \dotsc, \psi_m$ and modified psi classes
$\overline{\psi}_{m+1}, \dotsc, \overline{\psi}_n$ are thus our
gravitational descendants for the weight data $(1^m,0^{n-m})$.

This construction is extended to the case of $\beta=0, g=0$ in
\cite{Kock}. Here the modified psi classes are constructed with the
following twist: Start with $\oM_{0,n+3}(V,0)$ with three additional
distinguished marks $q_1,q_2,q_3$. For each of the other marks
$p_i,i\leq n$, define
\[
    \hat\pi_i : \oM_{0,n+3}(V,0) \rarr
    \oM_{0,\{q_1,q_2,q_3,p_i\}}\simeq\bP^1
\]
to be the map which forgets the sections not in
$\{q_1,q_2,q_3,p_i\}$ as well as the map to $V$. For $i\leq n$, the
definition is extended to this case to be
\[
\overline{\psi}_i:=\hat\pi_i^* (-2h),
\]
i.e. the pull-back of the class of degree $-2$ on $\bP^1$. A short
calculation is needed to see the connection, and we make it in the
following lemma.

\begin{lemma}
Consider the four pointed space $\oM_{0,A}$ with labeled points
$\{q_1,q_2,q_3,p_i\}$ and the weight data
$A=((1-\epsilon)^3,\epsilon)$ with $\epsilon<1/n$. Then
$\oM_{0,A}\simeq\bP^1$ and $\psi_{p_i,A}$ has degree $-2$.
\end{lemma}
\begin{proof}
We clearly have that $\oM_{0,A}\simeq\bP^1$ for any admissible $A$.
So we need only compute the degree of each psi class. The unweighted
basic pullback relation recalled in $\S$\ref{section:pullforget}
easily gives that each psi class on the unweighted space $\oM_{0,4}$
has degree 1. We note that $\Delta_{A}$ contains precisely the faces
for each vertex as well as the edges connecting the vertex $p_i$ to
each of the vertices $q_j$. Moreover, each divisor $D_{q_j,p_i}$ is
simply a point of $\bP^1$ and has degree 1. So application of
Corollary \ref{corollary:pullpsi} gives that
\[
\psi_{p_i,(1^4)}=\rho^*(\psi_{p_i,A})+D_{q_1,p_i}+D_{q_2,p_i}+D_{q_3,p_i},
\]
which we may easily solve to see $\psi_{p_i,A}$ has degree $-2$.
\end{proof}


\end{subsection}

\subsection{Generating functions}

\begin{definition}
  For fixed $n,V,\beta,\gamma_i$, we define the
  generating polynomial for the descendants to be
\[
e_{g,A}(\bt):=
\sum_{k_1,\dots,k_n}\des{\tau_{k_1}(\gamma_1)\cdots\tau_{k_n}(\gamma_n)}
{V,\beta}{g,A}\,\bt^\bk,
\]
and the exponential generating polynomial to be
\[
E_{g,A}(\bt):=
\sum_{k_1,\dots,k_n}
\frac{1}{\bk!}\des{\tau_{k_1}(\gamma_1)\cdots\tau_{k_n}(\gamma_n)}
{V,\beta}{g,A} \,\bt^\bk,
\]
and use the customary multi-index conventions that
$\bt:=(t_1,\dots,t_n)$, $\bk!:=k_1!\cdots k_n!$ and
$\bt^{\bk}:=t_1^{k_1}\cdots t_n^{k_n}$.
\end{definition}

These are polynomials because $\sum k_i$ is a constant based on the
choice of $n,V,\beta,\gamma_i$ by the dimension constraint on these
descendants, and hence there are only finitely many monomials in
these sums.

\begin{example}\label{example:genfunc}\quad
\begin{enumerate}
\item In the case of genus zero and weight $(1^n)$, it is well known that
\[
e_{0,n}=(t_1+\cdots+t_n)^{n-3}
\]

\item \label{example:permuto} It follows immediately from Lemma
\ref{lemma:permuto} and the above example, that
\begin{equation*}
e_{\Delta_{\cL_{n-2}}}=(t_1+t_2)^{n-3}
\end{equation*}
\item \label{example:projgen} It also follows immediately from Lemma
\ref{lemma:projpsi} that,
\[
E_{\Delta_{\bP^{r-2}}}=\frac{(t_1-t_2-\cdots-t_{r+1})^{r-2}}{(r-2)!}.
\]
\end{enumerate}
\end{example}

\section{Wall crossing formula}
\label{sec:wall-crossing-formula}

In this section, we consider the simple wall crossing given by
adding $\sigma$ to $\Delta_{\Ap}$ to form $\Delta_{A}$. As before,
denote the divisor corresponding to the crossing as $D_{\sigma}$,
and the reduction morphism $\rho_{A,\Ap}$ by $\rho$.

As in Theorem~\ref{thm:blowup}, we use an alternate notation $I$ for
$\sigma$, and $J$ for its complement.  The image of the divisor
$\Ds$ is given by the complex $\Delta_{A_J}=\Delta_{A_{\sigma}}$
which is attained from that of $\Delta_{\Ap}$ (or of $\Delta_{A}$)
by contracting the vertices in $\sigma$ to a disconnected vertex
which we label by $\sigma$. We also make the following notational
definition which we will use often.

\begin{definition} Let $\sigma$ be a collection of vertices of
  $\Delta_A$ and assume we are given a class $\tau_{k_i}(\gamma_i)$
  for each   vertex. Then we define
  $\boldsymbol{k_\sigma}:=\sum\limits_{i\in\sigma} k_i-\dim\sigma$ and
  $\boldsymbol{\gamma_\sigma}:=\prod\limits_{i\in\sigma}\gamma_{i}$.
\end{definition}


\begin{theorem}[Simple-wall Crossing]\label{theorem:wallcross}
We have the following wall crossing formula:

\begin{align*}
\des{\prod_{i=1}^{n}\tau_{k_i}(\gamma_i)}{V,\beta}{g,\Ap}&=\des{\prod_{i=1}^{n}\tau_{k_i}(\gamma_i)}{V,\beta}{g,\Am}+\\
&\qquad+(-1)^{\dim\sigma+1}\des{\tau_{k_{\sigma}}(\gamma_{\sigma})\prod_{j\not
\in\sigma}\tau_{k_j}(\gamma_j)}{V,\beta}{g,A_J}
\end{align*}
\end{theorem}

Before giving the proof of this result, we state some corollaries
and prove some preliminary results.
\begin{corollary}\label{corollary:gwinv}
The Gromov-Witten Invariants of $\oM_{g,A}(V,\beta)$ are equal to
the unweighted invariants for any weight $A$.
\end{corollary}
\begin{proof}
Indeed, in this case all $k_i=0$. Thus, $k_{\sigma}=-\dim\sigma<0$ and
therefore $\tau_{k_{\sigma}}=0$ by definition~\ref{defn:descendants}.
So the difference in
the wall-crossing formula is zero.
\end{proof}

\begin{corollary}
The generating polynomials are related by
\[
E_{g,\Ap}(\bt)-E_{g,A}(\bt)=(-1)^{\dim\sigma+1}
\Big(
\underbrace{\int\cdots\int }_{\dim\sigma}
E_{g,A_J}(\bt)
\underbrace{dt_{\sigma} \cdots dt_{\sigma}}_{\dim\sigma}
\Big)_{t_{\sigma}=\sum_{i\in \sigma}t_i}
\]
Here, the integration is homogeneous w.r.t. $dt_{\sigma}$, and hence
is uniquely defined, and after the integration we plug in
$t_{\sigma}=\sum_{i\in \sigma}t_i$.
\end{corollary}
\begin{proof}
  From formula \ref{theorem:wallcross}, this is a simple exercise on
  multinomial coefficients.
\end{proof}

\begin{lemma}\label{lemma:normBJ}
The normal bundle of ${\oM_{g,A_J}(V,\beta)}$ in
$\oM_{g,A}(V,\beta)$ is given by $N=\cO(-\psi_\sigma)^{\oplus
\dim\sigma}$.
\end{lemma}

Note that Theorem \ref{thm:blowup} implies that the normal bundle is
isomorphic to $E^{\oplus\dim\sigma}$ for \emph{some} line bundle
$E$. This lemma tells us what $E$ is.

\begin{proof}
  Let ${\oM}$ be the complete intersection of two sections $s_1$ and
  $s_2$ in the universal family $\cC_{\oM}$. Then $\oM$ is also the complete
  intersection of $s_1$ and $\cC_{\oM}.$ This implies that the normal
  bundle of ${\oM}$ in $s_2$ is isomorphic to the normal bundle of
  ${\oM}$ in $\cC_{\oM}$, which is $\cO_{\oM}(-\psi_1)$ by the definition of
  psi classes. On the locus where all the points $s_i$ coincide, we have
  $\cO(-\psi_1)=\cO(-\psi_{\sigma})$.

  Since $\oM_{g,A_J}(V,\beta)$ is a complete intersection of $r-1$ Cartier
  divisors, its normal bundle is a direct sum of $(r-1)$ of these line
  bundles.
\end{proof}

The general reference for the intersection theory on schemes that we
use is Fulton's \cite{Fulton_IT}. In particular, we recall from
\cite[Sec.6.2]{Fulton_IT} that for every regular embedding $i:X\to Y$ of
codimension $d$, morphism $f:Y'\to Y$, and a fiber diagram
\begin{equation*}
\xymatrix{
  X' \ar[d] \ar[r]^{j} & Y'\ar[d] \\
  X \ar[r]^i & Y
}
\end{equation*}
one has a refined Gysin homomorphism $i^!: A_k(Y')\to A_{k-d}(X')$
whose basic properties are listed in \cite[Thm.6.2]{Fulton_IT}. We
apply these properties to Artin stacks with finite stabilizers
instead of schemes. The Chow groups and intersection theory were
extended to the case of Artin stacks with affine stabilizer (this
includes our case of finite stabilizer) in the thesis of A. Kresch
\cite{Kresch}.

In the following lemma we use the notations of Definition~\ref{defn:IJ}.
We recall that we have a regular embedding of codimension one
$\mu:\oM_{0,A^+_I}\times\oM_{g,A_J}(V,\beta)\to\oM_{g,A^+}(V,\beta)$ and
that
$\oM_{0,A^+_I}=\bP^{r-2}$.

\begin{lemma}[Splitting Lemma]\label{lemma:shriek}
One has
$$
\mu^{!}[\oM_{g,A^+}(V,\beta)]\virt=
[\oM_{0,A_I^+}]\times[\oM_{g,A_J}(V,\beta)]\virt,
$$
\end{lemma}

This property does not obviously follow from the axioms of the virtual
fundamental cycle listed and checked in \cite{BayerManin}, so we
cannot just refer to the latter paper. In
particular, the divisor in \cite[6.3 Axiom 4c]{BayerManin} is bigger than
$D_{\sigma}$, since it contains all boundary components with
$\beta_1+\beta_2=\beta$, whereas $D_{\sigma}$ is only the component
corresponding to $\beta_1=0$, $\beta_2=\beta$.

\begin{proof}
Let $\Ds:=\oM_{0,A^+_I}\times\oM_{g,A_J}(V,\beta)$. Form the fiber
product diagram
\begin{equation*}
\xymatrix{
\Ds^1 \ar[d]_{q}\ar[rr]^{\mu_1} && \oM_{g,n}(V,\beta) \ar[d]^{\rho}\\
\Ds\ar[rr]^{\mu} && \oM_{g,A^+}(V,\beta)
}
\end{equation*}
with $\rho$ the reduction and $q$ the projection. Then
$\Ds^1=\oM_{0,|I|+1} \times \oM_{g,|J|+1}(V,\beta)$ is the divisor on
$\oM_{g,n}(V,\beta)$ corresponding to the graph with two vertices with
weights $(g_1,\beta_1)=(0,0)$ and $(g_2,\beta_2)=(g,\beta)$ and one
edge.  Since both $\mu$ and $\mu_1$ are regular embeddings of the same
codimension one, one has
\begin{displaymath}
  \mu_1^! [\oM_{g,n}(V,\beta)]\virt = \mu^! [\oM_{g,n}(V,\beta)]\virt
\end{displaymath}
by the compatibility property of Gysin morphisms, cf. \cite[Thm.6.2(c)]{Fulton_IT}.
  Therefore,
\begin{displaymath}
     \mu^! [\oM_{g,A^+}(V,\beta)]\virt =
     \mu^! \rho_* [\oM_{g,n}(V,\beta)]\virt =
     q_* \mu^!  [\oM_{g,n}(V,\beta)]\virt =
     q_* [\Ds^1]\virt.
\end{displaymath}
If $[\Ds^1]\virt = [\oM_{0,|I|+1}] \times [\oM_{g,|J|+1}(V,\beta)]\virt$, then
the latter equals $[\oM_{0,A^+_I}]\times[\oM_{g,A_J}(V,\beta)]\virt$, as
required. Therefore, the formula follows from the weight-1 case.

In the weight-1 case, this formula is well-known and follows by
applying \cite[Prop.7.2]{BehrendFantechi}, cf. \cite{Behrend}, Lemma
10 and proof of Axiom V (Isogenies).
\end{proof}

\begin{lemma}\label{lem:segre}
  One has
  \begin{eqnarray*}
    \rho_* \sum_{p>0} (-1)^{p-1}
    \left(D_{\sigma}^p \cap [\oM_{g,\Ap}(V,\beta)]\virt\right)
    \ =\  s(N) \cap [\oM_{g,A_J}(V,\beta)]\virt  \\
    =\ (1-\psi_{\sigma})^{-\dim\sigma} \cap [\oM_{g,A_J}(V,\beta)]\virt,
  \end{eqnarray*}
  where   $s(N)$ is the Segre class of the normal bundle.
\end{lemma}
\begin{proof}
  The latter equality is by the previous lemma:
  \begin{displaymath}
    s(N) = c\big( \cO(-\psi_{\sigma})^{\oplus \dim\sigma} \big)\inv
    = (1-\psi_{\sigma})^{-\dim\sigma}
  \end{displaymath}

  In the case when all
  $[\oM_{g,A}(V,\beta)]\virt=[\oM_{g,A}(V,\beta)]$ (for example if
  $V=\{pt\}$, or $g=0$ and $V$ is a convex variety), the former
  equality represents a well-known identity for the blowup of a
  complete intersection, see e.g.  \cite[3.3.4,6.7]{Fulton_IT}.  The
  formula applies to the blowup of any complete intersection in a
  scheme locally embeddable into a regular scheme. Since both sides
  behave functorially with respect to smooth covers, it holds for
  Artin stacks as well.


Here is the proof in the general case.
Recall that
  \begin{displaymath}
    D_{\sigma} = \oM_{0,A_I^+}\times \oM_{g,A_J}(V,\beta)
             = \bP^{r-2} \times \oM_{g,A_J}(V,\beta)
  \end{displaymath}
  is the exceptional divisor of a blowup. Let $\mu:D_{\sigma}\to
  \oM_{g,\Ap}(V,\beta)$ be the embedding and $\eta:D_{\sigma}\to \oM_{g,A_J}(V,\beta)$ be the
  projection. Then
  \begin{eqnarray*}
&&  \sum_{p>0} (-1)^{p-1} \rho_* \left(D_{\sigma}^p \cap [\oM_{g,\Ap}(V,\beta)]\virt \right) \\
&=& \sum_{p>0} (-1)^{p-1} \eta_* \left(D_{\sigma}^{p-1}\cap \mu^! [\oM_{g,\Ap}(V,\beta)]\virt \right)\\
&=& \sum_{p>0} \eta_* \left(c_1\big{(}\cO_{\bP(N)}(1)\big{)}^{p-1}
  \cap \mu^! [\oM_{g,\Ap}(V,\beta)]\virt\right) \\
&&\qquad\qquad\qquad\qquad\qquad\qquad
\text{because } \cO_{D_{\sigma}}(-D_{\sigma}) = \cO_{\bP(N)}(1) \\
&=& \sum_{q\ge0} \eta_*\left(c_1\big{(}\cO_{\bP(N)}(1)\big{)}^{q} \cap [\bP^{r-2}]\times[\oM_{g,A_J}(V,\beta)]\virt\right)\quad\text{ by Lemma \ref{lemma:shriek}}\\
&=& s(N) \cap [\oM_{g,A_J}(V,\beta)]\virt\quad \text{by the definition of the Segre classes.}
  \end{eqnarray*}
\end{proof}

\begin{proof}[Proof of Theorem \ref{theorem:wallcross}.]
We recall a few facts in preparing to apply the projection formula.
Theorem \ref{thm:blowup} gives that the reduction morphism $\rho$ is
given by a simple blowup along $\oM_{g,A_J}(V,\beta)$. It is
evident that the evaluation morphisms commute with reductions, and
so the projection formula allows us to push them forward unchanged.
In addition, Theorem \ref{thm:pullpsi} dictates that for
$i\in\sigma$, $\psi_{i,\Ap}=\rho^*(\psi_{i,A})+D_{\sigma}$, but the
remaining classes are pullbacks. We also note that for $i\in\sigma$,
whenever we restrict $\psi_{i,A}$ to $\oM_{g,A_J}(V,\beta)$ we get
$\psi_{\sigma}$.
We are now ready to compute.

\begin{eqnarray*}
&&
\des{\prod_{i\in\sigma}\tau_{k_i}(\gamma_i)  \prod_{j\not\in\sigma}\tau_{k_j}(\gamma_j)
}{V,\beta}{g,\Ap} = \\
&&
\prod_{i\in\sigma}(\rho^*\psi_{i,A}+D_{\sigma})^{k_i}\nu^*_{i,\Ap}(\gamma_i)
\prod_{j\not\in\sigma}(\rho^*\psi_{j,A})^{k_j}\nu^*_{i,\Ap}(\gamma_i)
\cap [\oM_{g,\Ap}(V,\beta)]\virt.
\end{eqnarray*}

Let us expand each $(\rho^*\psi_{i,A}+D_{\sigma})^{k_i}$ and
multiply them out. The only term that does not contain a positive
power of $D_{\sigma}$ is, by the projection formula and Definition
\ref{defn:virtual-class}, the descendant computed on $\oM_{g,A}(V,\beta)$. Let us
deal with the rest.

Note that all the terms with the evaluation classes and
with $\rho^*\psi_{j,A}$ are pullbacks from $\oM_{g,A}(V,\beta)$. Let us call this
part $\rho^*\tau_J(\gamma)$.

Now look at the remaining part, with $D_{\sigma}$ and
$\rho^*\psi_{i,A}$.  We now observe that, up to the sign $(-1)^{\sum
k_i+1}$, it is the homogeneous degree $\sum k_i$ part of
\begin{displaymath}
\prod_{i\in\sigma}\rho^*(1-\psi_{i,A})^{k_i}
\times
     \sum_{p>0} (-1)^p D_{\sigma}^{p-1}
\quad\text{applied to }  [\oM_{g,\Ap}(V,\beta)]\virt.
\end{displaymath}

By the projection formula $\rho_*(\rho^*\alpha\cap\beta)=
\alpha\cap\rho_*\beta$ and the previous lemma, we are reduced to
computing,
up to the sign $(-1)^{\sum k_i+1}$,
 the homogeneous degree $\sum k_i-\dim\sigma=k_{\sigma}$ part of
\begin{displaymath}
\prod_{i\in\sigma}(1-\psi_{i,A})^{k_i}
\times
     (1-\psi_{\sigma})^{-\dim\sigma}
\quad\text{applied to } \tau_J(\gamma)\cap [\oM_{g,A_J}(V,\beta)]\virt.
\end{displaymath}

But each
$\psi_{i,A}$ restrict to $\psi_{\sigma}$ on $\oM_{g,A_J}(V,\beta)$, so we need to
compute the degree $k_{\sigma}$ part of
\begin{displaymath}
\prod_{i\in\sigma}(1-\psi_{\sigma})^{k_i}
\times
     (1-\psi_{\sigma})^{-\dim\sigma}  = (1-\psi_{\sigma})^{k_{\sigma}},
\end{displaymath}
which is $(-1)^{k_{\sigma}}\psi^{k_{\sigma}}$. This gives the
formula. Note also that when $k_{\sigma}<0$, we get zero; all
monomials have nonnegative degree.
\end{proof}

\begin{definition}
  Let $\Sigma=\{\sigma\}$ be a partition of $\{1,\dotsc,n\}$ into
  a disjoint union of subsets. We say that $\Sigma$ is
  {\bf $\Delta_{A}$-admissible} if each $\sigma$ is in
  $\Delta_{A}.$ For each $\sigma\in\Sigma,$ we define
$\dim\Sigma:=\sum_{\sigma\in\Sigma}\dim\sigma$ and denote the number
of sets in the partition as $|\Sigma|.$ We denote the set of
$\Delta_{A}$-admissible
  partitions by $\Sigma(A).$ In addition, we notate $\Sigma(A,B)$ to be the set
  of partitions which are $\Delta_{B}$-admissible, but not
  $\Delta_{A}$-admissible.
\end{definition}

 We are now
ready to state a cousin of the string equation which reduces the
calculation of weighted descendants to that of the standard
unweighted descendants.

\begin{theorem}[Reduction to Unweighted Descendants]\label{theorem:absstring}
For any admissible weight data $A,$
\[
 \des{\prod_{i=1}^{n}\tau_{k_i}(\gamma_i)}{V,\beta}{g,A} =
 \sum_{\Sigma\in\Sigma(A)} (-1)^{\dim \Sigma}
  \des{\prod_{\sigma\in\Sigma}\tau_{k_\sigma}(\gamma_{\sigma})}{V,\beta}{g,|\Sigma|}
\]
\end{theorem}

\begin{proof}
  Pick a simplex $\sigma$ in $\Delta_A$ and apply formula
  (\ref{theorem:wallcross}) to get two complexes: one without
  $\sigma$, and one with $\sigma$ collapsed: it has one vertex instead
  of $\sigma$, disjoint from the rest.

  Now continue this inductively. The end result is the
  alternating sum over partitions of descendants on complexes which are
  disjoint unions of vertices, i.e. the unweighted descendants.
\end{proof}

\begin{corollary}
  The products of Miller-Morita-Mumford classes are expressed in the
  following way through the products of psi classes:
\[
\des{\kappa_{k_1-1}\cdots\kappa_{k_n-1}}{}{g,n}=
 \sum_{\text{\rm all partitions }\Sigma}
 (-1)^{\dim \Sigma}
  \des{\prod_{\sigma\in\Sigma}\tau_{k_\sigma} } {}{g,|\Sigma|}
\]
\end{corollary}
The inverse of this relation, ie expressing the psi numbers in terms
of the kappa numbers, is due to C. Faber and can be found in
\cite[1.13]{ArbarelloCornalba}.

\begin{corollary}\label{corollary:gendesc} For any $A\geq B$, we have:
\begin{align*}
  \des{\prod_{i=1}^{n}\tau_{k_i}(\gamma_i)}{}{g,B}&=\des{\prod_{i=1}^{n}\tau_{k_i}(\gamma_i)}{}{g,A}+\\
  &\qquad + \sum_{\Sigma\in\Sigma(A,B)} (-1)^{\dim\Sigma}
  \des{\prod_{\sigma\in\Sigma}\tau_{k_\sigma}(\gamma_{\sigma})}{}{g,|{\Sigma}|}
\end{align*}

\end{corollary}


\section{$A$-dilaton, $A$-string, $A$-divisor equations}
\label{sec:dilaton-string}

We are now in a position to use the results of
$\S$\ref{section:pullforget} to derive analogues of the well known
dilaton, string and divisor equations. For the remainder of this
section, we define $A$ to be of length $n+1$ and
$A':=A\setminus\{a_{n+1}\}$. We derive each first in the case
corresponding to the universal curve, and then in the case of a
symmetric weight corresponding to the complex $\Delta_{n+1,r}$ as
defined in Example \ref{example:compexamples}.\ref{example:symdef}.

The only property of the virtual fundamental cycle we need for these
computations is the property $[\cC_{g,A}(V,\beta)]\virt =
\pi_{A}^*[\oM_{g,A}(V,\beta)]\virt$ for the universal family, which
follows at once from the ``Forgetting Tails" property 
\cite[7.5(4)]{BehrendManin} and
Definition~\ref{defn:virtual-class} (cf. \cite[4.1,6.3]{BayerManin}).

\begin{subsection}{$A$-Dilaton Equation}

Recall that the unweighted Dilaton Equation states that
\[
\des{\tau_{1}\prod\limits_{i=1}^{n}\tau_{k_i}(\gamma_i)}{V,\beta}{g,n+1}=(2g-2+n)\des{\prod\limits_{i=1}^{n}\tau_{k_i}(\gamma_i)}{V,\beta}{g,n}
\]
The customary proof of this equality is to apply a push-pull type
argument with the forgetful map, using the pullback relation, and
reduce this to calculating the degree of a fiber;
we use it here as well.

We first give a version of the dilaton equation for the case
whenever the forgetful map corresponds to the map of the universal
curve. See Example
\ref{example:compexamples}.\ref{example:univcomplex}.

\begin{theorem}[Cone Dilaton Equation]
Assume $k_{n+1}=1$ and $\Delta_A=\Cone(\Delta_{A'})$. Then,
\begin{align*}
\des{\tau_{1}\prod_{i=1}^{n}\tau_{k_i}(\gamma_i)}{V,\beta}{g,A}&=
(2g-2)\des{\prod_{i=1}^{n}\tau_{k_i}(\gamma_i)}{V,\beta}{g,A'}
\end{align*}
\end{theorem}

\begin{proof}
Indeed, the pullback relation of Theorem~\ref{thm:forback} gives
that each psi class $1\leq i\leq n$ is given by the pullback. The
result thus follows from the projection formula and the degree of
fiber being $2g-2$.
\end{proof}

\begin{theorem}[Symmetric A-Dilaton
Equation]\label{theorem:symdilaton} Assume $k_{n+1}=1$ and $A$
corresponds to $\Delta_{n+1,r}$.  Then,
\begin{align*}
\des{\tau_{1}\prod_{i=1}^{n}\tau_{k_i}(\gamma_i)}{V,\beta}{g,A}&=
(2g-2)\des{\prod_{i=1}^{n}\tau_{k_i}(\gamma_i)}{V,\beta}{g,A'}+\\
&\quad\quad+\sum_{\sigma\in F(A,A')}(-1)^{\dim\sigma+1}\des{\tau_{{k_{\sigma}}}(\gamma_{\sigma})\prod_{j\not\in\sigma}\tau_{k_j}(\gamma_j)}{V,\beta}{g,A_{\sigma}}\\
\end{align*}
and $\Delta_{{A}_{\sigma}}$ is obtained from $\Delta_{n+1,r}$ by
combining the vertices in $\sigma$ (which includes $n+1$) to an
isolated vertex which is labeled by $\sigma$.
\end{theorem}

\begin{proof}
We make a simple reduction to the case of the universal curve which
is given above. The process of reducing the complex $\Delta_A$ to
$Cone(\Delta_{A'})$ is given by a reduction in bijection with
$F(A,A')$, by definition.  The only detail left to check is that the
complexes $\Delta_{A_\sigma}$ are independent of the order in which
these reductions are completed. We note by the definition of these
complexes as given in Example
\ref{example:compexamples}.\ref{example:sigmacompdef} that this is
equivalent to noting that for $\sigma_1\neq\sigma_2\in F(A,A')$, we
have that $\sigma_1\cap\sigma_2\neq\sigma_i$, ie there is no
containment among the $\sigma$'s. This is clear from the description
of $\Delta_{n+1,r}$.

Thus the result follows from applying the Wall-crossing formula
inductively to
\[
\des{\tau_{1}\prod_{i=1}^{n}\tau_{k_i}(\gamma_i)}{V,\beta}{g,A}
\]
and noting the first term in the statement follows from the case of
the universal curve which is given above.
\end{proof}

\begin{remark}
We note that the only ``symmetry" which is used in the above proof
is that there is no containment among the elements of $F(A,A')$, and
so the statement and proof are valid with these weaker hypotheses as
well. In the case where there are some containments among the
elements in $F(A,A')$, one must adjust this sum to account for the
extra faces (which will all contain the vertex labeled $\sigma$)
which appear in $A_{\sigma}$.
\end{remark}

\end{subsection}

\begin{subsection}{$A$-String Equation}
The unweighted string equation states
\[
\des{\tau_{0}\prod_{ i=1}^{n}
\tau_{k_i}(\gamma_i)}{V,\beta}{g,n+1}=\sum_{\ell=1}^{n}\des{\tau_{k_{\ell}-1}(\gamma_{\ell})\prod_{
i\neq \ell} \tau_{k_i}(\gamma_i)}{V,\beta}{g,n}
\]
We again give statements in the weighted case.

\begin{theorem}[Cone String Equation]
Assume $k_{n+1}=0$ and $\Delta_A=\Cone(\Delta_{A'})$. Then,
\[
\des{\tau_{0}\prod_{ i=1}^{n}\tau_{k_i}(\gamma_{i})}{V,\beta}{g,A}=0
\]
\end{theorem}
\begin{proof}
Indeed, this product pulls back from a space of a lower dimension
and is thus zero.
\end{proof}
\begin{theorem}[Symmetric $A$-String Equation]\label{proposition:symstring} Assume $k_{n+1}=0$ and $A$
corresponds to $\Delta_{n+1,r}$.  Then,
\[
\des{\tau_{0}\prod_{
i=1}^{n}\tau_{k_i}(\gamma_{i})}{V,\beta}{g,A}=\sum_{\sigma\in
F(A,A')}(-1)^{\dim\sigma+1}\des{\tau_{k_{\sigma}}(\gamma_{\sigma})\prod_{
i\not\in \sigma} \tau_{k_i}(\gamma_{i})}{V,\beta}{g,A_\sigma}
\]
and $\Delta_{{A}_{\sigma}}$ is obtained from $\Delta_{n+1,r}$ by
combining the vertices in $\sigma$ (which includes $n+1$) to an
isolated vertex which is labeled by $\sigma$.
\end{theorem}

\begin{proof}
The proof is virtually identical to that of Theorem
\ref{theorem:symdilaton} after noting the difference which appears
in the case of the universal curve.
\end{proof}

\end{subsection}

\begin{subsection}{$A$-Divisor Equation}
In the same spirit as the string and dilaton equations, there is the
usual unweighted divisor equation which states that for $D\in
A^1(V,\bQ)$
\begin{align*}
\des{\tau_{0}(D)\prod_{ i=1}^{n}
\tau_{k_i}(\gamma_i)}{V,\beta}{g,n+1}&=\int_{\beta}
D\cdot\des{\prod_{ i=1}^{n}
\tau_{k_i}(\gamma_i)}{V,\beta}{g,n}+\\
&\qquad\qquad+\sum_{\ell=1}^{n}\des{\tau_{k_{\ell}-1}(\gamma_{\ell}\cup
D)\prod_{ i\neq \ell} \tau_{k_i}(\gamma_i)}{V,\beta}{g,n}
\end{align*}

We state this now in the same cases as we stated the string and
dilaton equations. The proofs of each are very much in the same
spirit as previous proofs and are left to the reader.

\begin{theorem}[Cone $A$-Divisor Equation] Assume $k_{n+1}=0$ and $\Delta_A=\Cone(\Delta_{A'})$. For $D\in A^1(V,\bQ)$, we have
\begin{align*}
\des{\tau_{0}(D)\prod_{ i=1
}^{n}\tau_{k_i}(\gamma_{i})}{V,\beta}{g,A}&=\int_{\beta}
D\cdot\des{\prod_{ i=1}^{n}\tau_{k_i}(\gamma_{i})}{V,\beta}{g,A'}
\end{align*}
\end{theorem}
\begin{theorem}[Symmetric $A$-Divisor Equation]\label{proposition:symdivisor} Assume $k_{n+1}=0$ and $A$
corresponds to $\Delta_{n+1,r}$. For $D\in A^1(V,\bQ)$, we have
\begin{align*}
\des{\tau_{0}(D)\prod_{ i=1
}^{n}\tau_{k_i}(\gamma_{i})}{V,\beta}{g,A}&=\int_{\beta}
D\cdot\des{\prod_{ i=1}^{n}\tau_{k_i}(\gamma_{i})}{V,\beta}{g,A'}+\\
&\qquad+\sum_{\sigma\in
F(A,A')}(-1)^{\dim\sigma+1}\des{\tau_{k_{\sigma}}(\gamma_{\sigma}\cup
D)\prod_{ i\not\in \sigma}
\tau_{k_i}(\gamma_{i})}{V,\beta}{g,A_\sigma}
\end{align*}
and $\Delta_{{A}_{\sigma}}$ is obtained from $\Delta_{n+1,r}$ by
combining the vertices in $\sigma$ (which includes $j$) to an
isolated vertex which is labeled by $\sigma$.
\end{theorem}

\end{subsection}

\bibliographystyle{amsalpha}

\providecommand{\bysame}{\leavevmode\hbox to3em{\hrulefill}\thinspace}

\end{document}